\documentclass{cedram-aif}

\usepackage[francais, english]{babel}
\usepackage[utf8]{inputenc}
 \usepackage[T1]{fontenc}
 \usepackage{lmodern}

% ...

% ...

\title
[Low-lying zeros for Quaternion algebras]
{Low-lying zeros of L-functions for Quaternion Algebras}

\alttitle{Petits zéros de fonctions L pour les algèbres de quaternions}

\author{\firstname{Didier}  \lastname{Lesesvre}}

\address{School of mathematics (Zhuhai) \\
Zhuhai Campus, Sun Yat-Sen University \\
Tangjiawan, Zhuhai, 519082 \\
Guangdong, China (PRC)}

\email{lesesvre@math.cnrs.fr}

\thanks{I am infinitely indebted to my advisor, Farrell Brumley, for having trusted in me for handling this problem. I am grateful to the time Gergely Harcos and Philippe Michel granted me in carefully reading the thesis from which this article blossomed. I would like to thank Valentin Blomer, Andrew Corbett, Mikolaj Fraczyk, Élie Goudout and Guy Henniart for many enlightening discussions. I am also indebted to the referee for precisely reading the manuscript and providing comments and suggestions that contributed to improve the quality of this article. This work has been faithfully supported by the ANR 14-CE25 PerCoLaTor, the Fondation Sciences Mathématiques de Paris and the Deutscher Akademischer Austauschdienst.  At last, nothing would have been done without the warm and peaceful environment of the institutions which hosted me during these years: Université Paris 13, École Polytechnique Fédérale de Lausanne, Georg-August Universität and Sun Yat-Sen University.}

  
\altkeywords{représentation automorphe, fonction L, petits zéros, type de symétrie, conjecture de densité, opérateur de Hecke, algèbre de quaternions}


\usepackage{amsfonts, amssymb, amsmath} 
\usepackage{color}
\usepackage{graphicx}
\usepackage[all, cmtip]{xy}
\usepackage{eurosym}
\usepackage{amsthm}
\usepackage{stmaryrd}
\usepackage{appendix}
\usepackage{scrextend}
\usepackage{lipsum}
\usepackage[linkcolor=blue, citecolor=green, colorlinks=true]{hyperref}
\usepackage{array}
\usepackage{indentfirst}
\usepackage[Lenny]{fncychap}
\usepackage{thmtools}
\usepackage{thm-restate}
\usepackage{cleveref}

\newcolumntype{L}{>{$}c<{$}}

\newcommand{\R}{\mathbf{R}}
\newcommand{\C}{\mathbf{C}}

\newcommand{\N}{\mathbf{N}}
\newcommand{\A}{\mathbf{A}}

\newcommand{\Q}{\mathbf{Q}}

\newcommand{\GL}{\mathrm{GL}}

\newcommand{\PGL}{\mathrm{PGL}}

\renewcommand{\Re}{\mathrm{Re}}

\newcommand{\vol}{\mathrm{vol}}
\newcommand{\id}{\mathrm{id}}

\newcommand{\tr}{\mathrm{tr} \ }
\newcommand{\itex}{\item[$\bullet$]}
\newcommand{\tq}{\:: \:}
\renewcommand{\proof}{\textbf{Proof. }}

\newcommand{\p}{\mathfrak{p}}
\newcommand{\q}{\mathfrak{p}}
\newcommand{\dd}{\mathrm{d}}

\renewcommand{\det}{\mathrm{det}}

\renewcommand{\div}{\: | \:}
\renewcommand{\mod}{\mathrm{ mod } \ }

\newcommand{\ur}{\mathrm{ur}}
\renewcommand{\r}{\mathrm{r}}

\theoremstyle{plain}

\newcommand{\spec}{\mathrm{spec}}

\newcommand{\geom}{\mathrm{geom}}
\renewcommand{\tr}{\mathrm{tr }}
\renewcommand{\1}{\mathbf{1}}
\renewcommand{\q}{\mathfrak{q}}
\renewcommand{\d}{\mathfrak{d}}
\renewcommand{\c}{\mathfrak{c}}

\renewcommand{\p}{\mathfrak{p}}

\newcommand{\f}{\mathfrak{f}}
\newcommand{\n}{\mathfrak{n}}

\renewcommand{\dim}{\mathrm{dim } \ }
\renewcommand{\vol}{\mathrm{vol }}

\renewcommand{\GL}{\mathrm{GL}}

\renewcommand{\PGL}{\mathrm{PGL}}
\renewcommand{\O}{\mathcal{O}}

\renewcommand{\ell}{\mathrm{ell}}
\newcommand{\Pl}{\mathrm{Pl}}
\newcommand{\JL}{\mathrm{JL}}
\newcommand{\PPP}{\phi}

\renewcommand{\proof}{\textit{Proof. }}

\newcommand{\rk}{\noindent \textit{Remark. }}

\begin{document}

\begin{abstract}
The density conjecture of Katz and Sarnak predicts that, for natural families of L-functions, the distribution of zeros lying near the real axis is governed by a group of symmetry. In the case of the universal family of automorphic forms on a totally definite quaternion algebra,  we determine the associated distribution for a restricted class of test functions in the analytic conductor aspect. In particular it leads to non-trivial results on densities of non-vanishing at the central point. 
\end{abstract}

\begin{altabstract}
Pour des familles naturelles de fonctions L, La conjecture de densité de Katz et Sarnak prédit que la distribution des zéros proches de l'axe réel est régie par un groupe de symétrie. Dans le cas de la famille universelle des algèbres de quaternions totalement définies, nous déterminons la distribution associée pour une classe explicite de fonctions test, uniformément lorsque le conducteur analytique croît. En particulier, cela mène à des résultats non-triviaux sur les densités de non-annulation aux valeurs centrales.
\end{altabstract}

\maketitle

\section{Introduction}
\label{sec:intro}

\subsection{Statement of results}
\label{subsec:landscape}

Let $F$ be a number field of degree $d$ over $\Q$ and $\A$ the ring of adeles of $F$. We consider a division quaternion algebra $B$ over $F$, and write $R$ for the places of $F$ where $B$ is not split. We introduce the group of projective units $G = Z \backslash B^\times$, where $Z$ denotes the center of $B^\times$. Let $\mathcal{A}(G)$ denote the universal family of $G$, that is the set of all infinite dimensional irreducible automorphic representations of the group $G(\A)$. It embeds, via the Jacquet-Langlands correspondence, into the universal family of $\PGL(2)$ made of all its cuspidal automorphic representations. Following Sarnak \cite{sarnak_definition_2008}, a deep understanding of $\mathcal{A}(G)$ is of fundamental importance in the theory of automorphic forms. 

In order to make sense of problems on average for $\mathcal{A}(G)$, it is necessary to truncate the universal family to a finite set. We do so by bounding the analytic conductor $c(\pi)$ of Iwaniec and Sarnak \cite{iwaniec_perspectives_2000}. The truncated universal family may then be introduced as
\begin{equation}
\label{UF}
\mathcal{A}(Q) = \{\pi \in \mathcal{A}(G) \tq   c(\pi) \leqslant Q\}, \qquad Q \geqslant 1.
\end{equation}

In a previous work \cite{lesesvre_counting_2020}, we have shown that $\mathcal{A}(Q)$ is a finite set and determined arithmetic statistics on this universal family when $Q$ grows to infinity. Namely, an asymptotic expansion for its cardinality was given, as well as global and local equidistribution results. 

\begin{thm}[Counting law for quaternion algebras]\label{thm-count}
Introduce $\delta_F = 2({1+[F:\Q]})^{-1}$. There exists $C > 0$ such that, for any $Q \geqslant 1$, 
\begin{equation*}
|\mathcal{A}(Q)| =  C Q^2 + 
\left\{
\begin{array}{cl}
O\left(Q^{1+\varepsilon}\right) & \text{ if } F = \Q \text{ and $B$ totally definite, for all } \varepsilon >0; \\[.5em]
O\left(Q^{2-\delta_F} \right) & \text{ if $F \neq \Q$} \text{ and $B$ totally definite} ;\\[.3em]
O\left(\frac{Q^2}{\log Q} \right) & \text{ if $B$ is not totally definite}.
\end{array}
\right.
\end{equation*}
\label{thm1}
\label{deltaF}
\end{thm}

The exact form of the constant $C$ is given in \cite[Eq. (8)]{lesesvre_counting_2020} in a geometric meaningful form. The knowledge of this cardinality opens the path to other statistical results. By refining the methods leading to the above theorem, mainly based on the construction of a suitable test function for the Arthur-Selberg trace formula, we address here the question of the distribution of the low-lying zeros of the associated L-functions. We briefly state the two main results of this paper before presenting the more detailed setup in the next sections.

Let $\pi \in \mathcal{A}(G)$. Let $\phi$ be an even Schwartz function on $\R$ with compactly supported Fourier transform. In particular it admits analytic continuation to the whole the complex plane. The one-level density attached to $\pi$ is defined by the distribution
\begin{equation}
\label{density-llz}
D(\pi, \phi) = \sum_{\gamma_\pi} \phi\left(\tilde{\gamma}_\pi\right),
\end{equation}

\noindent where the sum runs over nontrivial zeros of the $L$-function $L(s, \pi)$, normalized as in Section \ref{sec:norm}.

\begin{restatable}{thm}{TS}
\label{thm:ts}
\label{thmD}
Let $B$ a totally definite quaternion algebra. For every even and Schwartz class function $\phi$ on $\R$ with Fourier transform compactly supported in $(-2/3, 2/3)$, we have
\begin{equation}
\label{ts}
\frac{1}{|\mathcal{A}(Q)|} \sum_{\pi \in \mathcal{A}(Q)} D(\pi, \phi) \xrightarrow[Q \to \infty]{} \int_\R \phi(x) W_O(x) \mathrm{d}x =  \widehat{\phi}(0) + \frac{1}{2} \phi(0).
\end{equation}

\noindent  The distribution density is therefore given by the function $W_O = 1 + \frac{1}{2}\delta_0$.  In particular, the type of symmetry of the inner forms of \, $\mathrm{PGL}(2)$ is the one of the orthogonal group.
\end{restatable}

Statistics on the distribution of low-lying zeros of $L$-functions are known to lead to results concerning vanishing at the central point, following the ideas of Iwaniec, Luo and Sarnak \cite{iwaniec_low_2000}. Introduce the proportion of automorphic representations with vanishing at the central point of order $m$, that is to say
\begin{equation}
p_m(Q) = \frac{1}{|\mathcal{A}(Q)|} \# \left\{ \pi \in \mathcal{A}(Q) \ : \ \underset{{s=1/2}}{\mathrm{ord}} \ L(s, \pi) = m \right\}, \qquad m \in \N.
\end{equation}

\begin{restatable}{coro}{nonvanishing}
\label{coro:coro}
Let $B$ a totally definite quaternion algebra and assume the generalized Riemann hypothesis. We have
\begin{equation}
\liminf_{Q \to \infty} \ \sum_{m \geqslant 1} mp_m(Q)  \leqslant 2.
\end{equation}
\end{restatable}

\subsection{One-level densities}

There is a long-lasting history about zeros of L-functions, that remain even today a challenging field of investigation. The lack of tools to grasp such zeros pointwise leads to search for results on average. The theory of random matrices \cite{mehta_random_2004} is a glass through which
understand the field of statistics on zeros of $L$-functions. Indeed, the zeros of families of $L$-functions behave strikingly like the eigenangles of classical groups of random matrices, and this heuristics serves as a guide for the $L$-function world.  For a wide insight, we refer to the general survey of Sarnak, Shin and Templier\cite{sarnak_families_2016}.
 
Rudnick and Sarnak \cite{rudnick_zeroes_1996} showed that the pair correlation of spacings between zeros of automorphic L-functions matches the pair correlation of spacings between eigenvalues of classical groups of random matrices.  The universality of the pair correlations is surprising for it is blind to the differences between the classical groups. A further disappointment with correlations is that they are unsensitive to finitely many modifications of the zeros, and in particular do not give any importance to zeros usually of arithmetic significance, typically the central point. Altogether, it can be expected that other statistics are able to distinguish between them.

\subsubsection{One-level density for matrices}

The correlation statistics take into account \textit{all} the eigenangles, since they consider only the distribution of spacings between them. Katz and Sarnak \cite{katz_random_1999} broke this universality, turning their interest
towards statistics concentrated on \textit{small} eingenangles. They proved that the average
density of these small eigenangles over a family differs depending on the group. Let $\phi$ be an even Schwartz function on $\R$ and and $A$ be a diagonalizable unitary matrix of one of the classical groups $G(N)$. The one-level density attached to $A \in G(N)$, where $G(N)$ is among $\mathrm{U}(N), \mathrm{Sp}(N), \mathrm{O}(N), \mathrm{SO}(2N)$ or $ \mathrm{SO}(2N+1)$, is
\begin{equation}
\label{density-llz}
D(A, \phi) = \sum_{\theta_A} \phi\left(\tilde{\theta}_A\right),
\end{equation}

\noindent where the sum runs over the eigenangles $\theta_A$ of $A$, and $\tilde{\theta}_A = \frac{N}{2\pi} \theta_A$.

\begin{thm}[Katz-Sarnak]
For the classical groups $G(N)$, for every even Shcwartz function $\phi$ on $\R$ whose Fourier transform is compactly supported, we have
\begin{equation}
\int_{G(N)} D(A, \phi) \mathrm{d}A \xrightarrow[N \to \infty]{} \int_{\R} W_G(x)\phi(x)\dd x,
\end{equation}

\noindent where $\mathrm{d}A$ is a normalized Haar measure on $G(N)$, and the densities functions $W_G$ on $\R$ are defined by
\begin{align*}
W_{\mathrm{U}}(x) & = 1 \\ 
W_{\mathrm{Sp}}(x) & = 1 - \frac{\sin 2 \pi x}{2 \pi x} \\ 
W_{\mathrm{SO(even)}}(x) & = 1 + \frac{\sin 2 \pi x}{2 \pi x} \\
W_{\mathrm{SO(odd)}}(x) & = 1 - \frac{\sin 2 \pi x}{2 \pi x} + \delta_0(x) \\
W_{\mathrm{O}}(x) & = \frac{1}{2} \left( W_{\mathrm{SO(even)}}(x) + W_{\mathrm{SO(odd)}}(x) \right) = 1 + \frac{1}{2} \delta_0(x) .
\end{align*}
\end{thm}

The
fact that the limit is no more universal but does depend on the group gives rise to the notion of \textit{type of symmetry} of a family of random matrices.

\subsubsection{One-level density for $L$-functions}
\label{sec:norm}

Following the enlightening analogy with random matrices, it can be expected
that the one-level density of the zeros attached to every reasonable family of $L$-functions behaves as the one-level density of the eigenangles of the classical groups of random matrices, in other words
that the behavior of low-lying zeros of $L$-functions is modeled by
the classical groups. 

Let $\pi$ be an automorphic representation with an associated notion of $L$-function $L(s, \pi)$. Consider its nontrivial zeros written in the form
\begin{equation}
\rho_\pi = \frac{1}{2}+i\gamma_\pi.
\end{equation}
Here the $\gamma_\pi$ are a priori complex numbers without assuming the Riemann
hypothesis.  Renormalize the mean spacing (in the case where the $\gamma_\pi$ are real) of the low-lying zeros to $1$ by setting
\begin{equation}
\label{zeros}
\tilde{\gamma}_\pi = \frac{\log c(\pi)}{2\pi} \gamma_\pi.
\end{equation}

Let $\phi$ be an even Schwartz function on $\R$ whose Fourier transform is compactly supported, in particular it admits an analytic continuation to all $\C$. Let $\pi$ be an automorphic representation. The one-level density attached to $\pi$ is defined by
\begin{equation}
\label{density-llz}
D(\pi, \phi) = \sum_{\gamma_\pi} \phi\left(\tilde{\gamma}_\pi\right).
\end{equation}

The first result estimating the one-level density is given by Özlük and Snyder \cite{ozluk_small_1993} in 1993 for $L$-functions attached to Dirichlet characters. Since then, a wide literature has been published concerning the statistical behavior of low-lying zeros of families of $L$-functions \cite{duenez_low_2006, goldfeld_gl3_2013, iwaniec_low_2000, liu_low-lying_2017, rubinstein_low-lying_2001}. This led Katz and Sarnak \cite{katz_zeroes_1999} to formulate the so-called \emph{density conjecture}, claiming the same universality for the types of symmetry of families of $L$-functions than those arising for classical groups of random matrices.
\begin{conj}[Density conjecture]
\index{density conjecture}
Let $\mathcal{F}$ be a family of automorphic representations in the sense of Sarnak and $\mathcal{F}_Q$ a finite truncation increasing to $\mathcal{F}$ when $Q$ grows. Then for all even Schwartz function on $\R$ with compactly supported Fourier transform, there is one classical group $G$ among $\mathrm{U}$, $\mathrm{SO(even)}$, $\mathrm{SO(odd)}$, $\mathrm{O}$ or $\mathrm{Sp}$ such that
\begin{equation}
\frac{1}{|\mathcal{F}_Q|} \xrightarrow[Q\to\infty]{} \int_\R \phi(x) W_G(x)\dd x.
\end{equation}
The family $\mathcal{F}$ is then said to have the type of symmetry of $G$.
\end{conj}

\rk  For families of $L$-functions associated to algebraic varieties over function fields, the type of symmetry is determined by the monodromy of the family \cite{katz_random_1999}, shedding light on the reason why zeros of $L$-functions are governed by groups of random matrices. No such analogue is known for number fields. 




\subsection{Type of symmetry of quaternion algebras}
\label{subsec:one-level}
\label{subsec:intro-ts}

Considering the statistics on low-lying zeros of $L$-functions attached to the universal family of quaternion algebras, the problem is to determine whether or not the averaged one-level density on $\mathcal{A}(Q)$ admits a limit and unveils a type of symmetry according to the density conjecture. The following statement answers positively to this question and partially determines the type of symmetry of quaternion algebras. 
\TS*


\noindent \textit{Remark.} It is natural to relate the above result to the two important papers~\cite{iwaniec_low_2000} and~\cite{shin_sato-tate_2016}.
\begin{itemize}
\item[(i)] The landmark work in the Katz-Sarnak philosophy is the one of Iwaniec, Luo and Sarnak \cite{iwaniec_low_2000} in the case of classical automorphic forms. By the Jacquet-Langlands correspondence, the universal family we consider, $\mathcal{A}(G)$, can be embedded in the universal family made of all the cuspidal automorphic representations of $\mathrm{GL}(2)$. When $B$ is the quaternion algebra over $\Q$ ramified at the prime place $q$ and at the archimedean place, the automorphic representations we consider are in one-to-one correspondence with classical holomorphic cusp forms of level $q$, so that the result of Iwaniec, Luo and Sarnak provides and instance of Theorem \ref{thmD}. Their result displays a better bound for the support allowed for the Fourier transform (viz. $2$ instead of $2/3$). This feature comes from the fact that their setting is restricted to the special case of holomorphic cusp forms, allowing to use the Petersson formula and specific bounds on Kloosterman sums. Theorem \ref{thmD} is much more general and not only encapsulates Hilbert modular forms on general number fields, but also Maass forms.
\item[(ii)]  Despite the wide generality of \cite{shin_sato-tate_2016}, Shin and Templier are not able to display a quantitatively explicit bound for the limiting Fourier support and do require automorphic representations to have prescribed ramification.  Theorem \ref{thmD} allows automorphic representations to be ramified at an arbitrarily large and not uniformly bounded number of places.
\end{itemize} 

An important caveat ought to be mentioned concerning the orthogonal types of symmetry. The density conjecture postulates results for Schwartz functions with arbitrary compactly supported Fourier transform. Assuming this conjecture, proving the convergence for a narrower class of allowed Fourier supports may determine uniquely the conjectural type of symmetry. However this is not the case for all the supports, and an uncertainty remains in the case or supports smaller than $(-1, 1)$. Indeed, the Plancherel formula yields
\begin{equation}
\int_\R \phi(x) W(x) \dd x = \int_\R \widehat{\phi}(x) \widehat{W}(x) \dd x.
\end{equation}

\noindent Looking at the Fourier transforms of the densities, introducing $\eta$ the characteristic function of $[-1, 1]$, direct computations leads to
\begin{align*}
\widehat{W}_{\mathrm{U}}(x) & = \delta_0(x) \\ 
\widehat{W}_{\mathrm{Sp}}(x) & = \delta_0(x) - \frac{1}{2}\eta(x) \\ 
\widehat{W}_{\mathrm{SO(even)}}(x) & = \delta_0(x) + \frac{1}{2}\eta(x) \\
\widehat{W}_{\mathrm{SO(odd)}}(x) & = \delta_0(x) -\frac{1}{2}\eta(x) + 1 \\
\widehat{W}_{\mathrm{O}}(x) & = \delta_0(x) + \frac{1}{2}.
\end{align*}

\noindent Unfortunately, we notice that the three orthogonal types of symmetry, viz. $\widehat{W}_{\mathrm{O}}$, $\widehat{W}_{\mathrm{SO(even)}}$ and $\widehat{W}_{\mathrm{SO(odd)}}$, are indistinguishable in $(-1, 1)$. Therefore, Theorem \ref{thmD} only partially determines the type of symmetry of the universal family of quaternion algebras.  However, following Miller \cite{Miller_1-and}, determining the two-level density of the low-lying zeros of $\mathcal{A}(G)$ is enough to determine the type of symmetry for an arbitrarily small support of the Fourier transform of the test function. This computation of the limiting two-level density can be handled explicitly, and essentially amounts to squaring the expression \eqref{old} and to carry on the strategy of the present paper. The formal computation following from there is carried out \textit{mutatis mutandis} in \cite[Section 4.3]{Miller_1-and}. It can be observed that the determination of the precise type of symmetry, among the three orthogonal ones, depends on the proportion of automorphic forms having sign $-1$ in the functional equation. Such information is not available in the level of generality of this paper; however, with mild assumption on the conductor, Martin \cite{martin_refined_2018} shows (despite an interesting bias) that the limiting proportion is $1/2$, pleading for an orthogonal type of symmetry for $\mathcal{A}(G)$. The family restricted to positive (resp. negative) signs in the functional equation has an even (resp. odd) orthogonal type of symmetry, see \cite[Theorem 3.2]{Miller_1-and}. This formal argument provides evidence for the conjectural type of symmetry; however a careful study of the spectral terms involved in the trace formula is needed to ensure actual convergence for non-trivial supports of the Fourier transform.

\noindent \textit{Remarks.} There are some directions in which this result generalizes.
\begin{itemize}
\item[(i)] The centerless assumption has been made for convenience and the same proof carries on to the general central character with minor modifications in Section \ref{subsubsec:case}. In this case, the type of symmetry is expected to be unitary.
\item[(ii)] The totally definite assumption is more significant, but the result is still expected to hold in the case of any quaternion algebra. Let us mention some of the extra difficulties that should arise. For a non totally definite division quaternion algebra, the archimedean spectrum is no more discrete and the selecting function has to be approximated by a smoothing procedure. This can be done only on the tempered part of the spectrum, by Paley-Wiener type theorems, so that the contribution of the archimedean complementary spectrum has to be shown to contribute as an error term. This has already been done by the author for counting and equidistribution results \cite{lesesvre_counting_2020} and the same method should be adaptable to extend Theorem \ref{thmD} to this setting, see the author's thesis \cite{lesesvre_thesis} for the details. 
\item[(iii)] For the case of $\GL(2)$, the automorphic spectrum admits a continuous part that should also be carefully taken care of. The methods Brumley and Mili\'{c}evi\'{c} \cite{brumley_counting_2016} developped  to prove the counting law for the universal family of $\GL(2)$ provides the necessary tools. However, these adaptations are far from trivial and require highly technical analysis of these extra spectral terms.
\item[(iv)] The result on the density of vanishing at the central point, in addition of verifying the whole density conjecture, are strong motivations to strengthen the bounds on the support of the Fourier transform in Theorem \ref{thmD}, that would yield better results. However, this would require stronger estimates on orbital integrals that are not available in the current literature.
\end{itemize}

\subsection{Organization of this article}

In Section \ref{sec:spectral-sums} we recall the basic properties of the automorphic L-functions and we reformulate the one-level density in terms of distributions of  prime numbers, via an explicit formula. The notations and measures are introduced and the analytic conductor is precisely defined there. The universal family can be decomposed into harmonic subfamilies  by fixing spectral data, easier to understand, as explained in Section \ref{sec:decomposition}. The high orders contributions in the explicit formula are shown to be negligible in Section \ref{sec:high-contrib}. The remaining terms are related to the Hecke eigenvalues in Section \ref{sec:traces-hecke}, leading to a reformulation of the low orders terms as a spectral side of the Selberg trace formula for a suitable test function, constructed in Section \ref{sec:stf}. A precise analysis of the geometric side leads to unveiling the nontrivial contribution of the low orders terms to the type of symmetry in Section \ref{sec:low-order}, achieving the proof of Theorem \ref{thmD}. Finally, the result on the density of non-vanishing at the central point is derived from the type of symmetry in Section \ref{sec:non-vanishing}.

\section{Groundwork}
\label{sec:spectral-sums}

We denote by $v$ the places of $F$, $\p$ the non-archimedian ones, and $\O_\p$ the ring of integers of $F_\p$ for a finite place $\p$. The finite set $R$ of ramification places of $B$ determines it up to isomorphism, and is assumed to contain all the archimedean places. From now on, Latin letters $q, d, m$, etc. will denote rational integers, while Gothic letters $\q, \d, \mathfrak{m}$, etc. will denote ideals of integer rings. 

\subsection{Analytic conductor}
\label{subsec:conductor}
\label{sec:conductors}

In order to make sense of the problem, we need to define precisely the notion of size we choose for representations. It is the analytic conductor, which we introduce in this section. We will work with $B^\times$ more than with~$G$, for it lightens notations. This local convention makes no harm, for we view a representation $\pi$ of $G(\A) = PB^{\times}(\A)$ as a representation of $B^\times(\A)$ with trivial central character. By Flath's theorem, an irreducible admissible representation of $B^\times(\A)$ decomposes in a unique way as a restricted tensor product $\pi = \otimes_v \pi_v$ of irreducible smooth representations where almost every component $\pi_v$ is unramified. We therefore want first to define the conductor for the local components $\pi_v$.

The Jacquet-Langlands correspondence allows to reduce to the $\GL(2)$ case, and in this one only infinite-dimensional representation arise. Indeed, since the universal family excludes global characters, a representation $\pi$ in it is generic. The Jacquet-Langlands correspondence preserves genericity hence, as shown on the diagram below, associates to $\pi$ a generic representation $\JL(\pi)$ of $\GL(2)$, thus also its local components $\JL(\pi)_v$. These local components are also the images by the local Jacquet-Langlands correspondence $\JL_v(\pi_v)$ of the local components of $\pi$. 
$$
\xymatrix{
\, \substack{\pi \in \mathcal{A}({B^\times})} \quad \ar[rr]^{\mathrm{JL}} \ar[d]_v & & \quad \, \substack{\mathrm{JL}(\pi) \in \mathcal{A}(\GL(2)) \\ \mathrm{generic}} \ar[d]^v \\
\pi_v \quad \ar[rr]^{\mathrm{JL}_v}_{\id \ \mathrm{if} \ v \notin R} & & \quad \substack{\mathrm{JL}(\pi)_v \\ \mathrm{generic}}
}
$$

At split places, the local Jacquet-Langlands correspondence is the identity, for then $B_\p^\times \simeq \GL(2, F_\p)$. The correspondence is unique, thus the local components at split places $\pi_v$ are generic hence infinite-dimensional, proving the claim. 

\subsubsection{Non-archimedian case}

For finite split places $\p$, by definition $B_{\p} \simeq M\left(2, F_{\p}\right)$ so that $B^\times_{\p}$ is isomorphic to $\GL(2, F_\p)$. The notion of local conductor for irreducible smooth infinite-dimensional representations of $\GL(2)$ has been introduced by Casselman \cite{casselman_results_1973}. Consider the sequence of compact open congruence subgroups, for $r \geqslant 0$, 
\begin{equation}
\label{filtration}
K_{0, \p}\left(\p^r\right) = 
\left\{
g \in \GL\left(2, \O_{\p}\right) \tq 
g \equiv
\left(
\begin{array}{cc}
\star & \star \\
0 & \star
\end{array}
\right)
\mod \p^r
\right\} \subseteq B^\times_{\p}.
\end{equation}

The multiplicative and analytic conductors of an  irreducible admissible infinite-dimensional representation $\pi_\p$ of $B^\times_\p$ with trivial central character are then respectively defined by
\begin{equation}
\label{conductor-XXX}
\c(\pi_\p) = \p^{\f(\pi_{\p})} \quad \text{and} \quad  c(\pi_\p) = N\c(\pi_\p),
\end{equation}

\noindent where 
\begin{equation}
 \f(\pi_{\p}) = \min \left\{r \in \N \tq \pi_{\p}^{K_{0, \p}\left(\p^r\right)} \neq 0\right\}.
\end{equation}

\noindent The existence of the conductor is guaranteed by the work of Casselman \cite{casselman_results_1973}, who also states that the growth of the dimensions of the fixed vector spaces are given by
\begin{equation}
\dim \pi_\p^{K_{0, \p}\left(\p^{\f(\pi_\p)+i}\right)} = i+1, \qquad i \geqslant 0.
\label{one}
\end{equation}

\subsubsection{Archimedian case}

The archimedian part of the conductor is introduced by Iwaniec and Sarnak \cite{iwaniec_perspectives_2000}. It is built on the archimedean factors completing the L-functions associated to automorphic representations. The archimedean L-factors are of the form
\begin{equation}
L(s, \pi_v) = \prod_{j=1}^2 \Gamma_v(s-\mu_{j,\pi}(v)),
\end{equation}

\noindent where $\Gamma_v(s) = \pi^{-s/2}\Gamma(s/2)$ and the $\mu_{j,\pi}(v)$ are complex numbers. The local analytic conductor is then locally defined to be, for $v \div \infty$, 
\begin{equation}
c_v(\pi) = \prod_{j=1}^2 \left(1 + |\mu_{j, \pi}(v)|\right).
\end{equation}

\subsubsection{Non-split case}

Via the Jacquet-Langlands correspondence, the non-split case is reduced to the already treated split one. The conductor of an irreducible admissible representation $\pi_{v}$ of $B^\times_{v}$ is defined as the conductor of its Jacquet-Langlands transfer
\begin{equation}
c\left(\pi_{v}\right) = c\left(\JL\left(\pi_{v}\right)\right).
\end{equation}

\subsubsection{Characters}

For now conductors have been defined only for generic representations. However, characters can arise as local components at ramified places as discussed above. Every character of $B_\p^\times$ is a composition
\begin{equation}
B_\p^\times \longrightarrow F_\p^\times \longrightarrow \C,
\end{equation}

\noindent where the first application is the reduced norm, and the second one a character of $F_\p^\times$. In other words, every character of $B_\p^\times$ is of the form $\chi_0 \circ N$ where $\chi_0$ is a character of $F_\p^\times$ and $N$ the reduced norm on $B_\p^\times$. In order to stay consistent, define the conductor of a local character at a ramified place as the conductor of its Jacquet-Langlands embedding in $\GL(2)$. Since the character $\chi_0 \circ N$ is sent to the twisted Steinberg representation $\mathrm{St} \otimes \chi_0$, it follows explicitly
\begin{equation}
\c(\chi_0 \circ N) = \left\{
\begin{array}{cl}
\p & \text{if $\chi_0$ unramified;} \\
\c(\chi_0)^2 & \text{if $\chi_0$ ramified.}
\end{array}
\right.
\end{equation}

\subsubsection{Global analytic conductor}

For an irreducible admissible representation of $B^\times(\A)$ decomposed into $\pi = \otimes_{v} \pi_{v}$, introduce its global analytic conductor
\begin{equation}
c\left(\pi\right) = \prod_{v} c\left(\pi_v\right).
\end{equation}

\noindent This gives a well-defined notion, for the $\pi_v$ are almost everywhere unramified, thus of conductor one. It extends to a definition for representations of $G(\A)$, viewed as automorphic representations of $B^\times(\A)$ with trivial central characters.

\subsection{Normalization of measures}
\label{subsec:measures}
\label{sec:plancherel-formulas}
\label{subsec:Haar-measures}

At the non-archimedean places, the measure taken on $G_\p$ is the Haar measure $\mu_\p$ normalized so that $K_\p = \PGL(2,\O_{\p})$, in the split case, or $K_\p = \mathfrak{o}_{\p}^\times$ the units of a maximal order of $B_\p$, in the non-split case, gets measure one. This normalization is independent \cite{hull} of the chosen maximal order. For the archimedean places, we choose the Haar measure normalized so that the maximal compact subgroup gets measure one. 

We now turn to the associated local dual groups. Denote $\mathcal{H}(G_v)$ the Hecke algebra of $G_v$, that is the algebra consisting of compactly supported complex-valued functions on $G_v$, locally constant at finite places, smooth at archimedian ones. Let $\mathcal{H}(G(\A))$ be the Hecke algebra of $G(\A)$. It is the algebra generated by the restricted products $\phi = \prod_v \phi_v$, where $\phi_v$ is a function of $\mathcal{H}(G_v)$ and almost every local component $\phi_\p$ is equal to $\1_{K_{\p}}$. For such a function $\phi \in \mathcal{H}(G(\A))$, we extend the action of $\pi$ to $\mathcal{H}(G(\A))$, letting $\pi(\phi)$ act by the mean action of $\pi$ weighted by $\phi$, that is to say
\begin{equation}
\pi(\phi) = \int_{G(\A)} \phi(g) \pi(g) \dd g.
\end{equation}

\noindent This defines a Hilbert-Schmidt integral operator of trace class, thus we can define its Fourier transform by
\begin{equation}
\widehat{\phi}(\pi) = \tr \ \pi(\phi) = \tr \left( v \mapsto \int_{G} \phi(g)\pi(g)v \dd g \right).
\end{equation} 

The unitary dual group $\widehat{G}_v$ is endowed with its usual Fell topology and Plancherel measure associated with the measure chosen on $G_v$: it is the unique positive Radon measure $\mu^{\Pl}_v$ on $\widehat{G}_v$ such that the Plancherel inversion formula of Harish-Chandra holds, \textit{i.e.} for functions $\phi_v$ in the Hecke algebra $\mathcal{H}(G_v)$, we have
\begin{equation}
\int_{\widehat{G}_v} \widehat{\phi}_v(\pi_v) \dd\mu_v^{\Pl}(\pi_v) = \phi_v(1).
\label{plancherel}
\end{equation}

\noindent From now on, integrals on $\widehat{G}_v$ will be written with the convention that $\dd\pi_v = \dd\mu_v^{\Pl}(\pi_v)$, leading to no ambiguity. On $\widehat{\Pi} = \prod_v \widehat{G}_v$ we consider the product topology and the Plancherel measure, denoted by $\mu^\Pl$ and given by the product of the local ones. 

\subsection{L-functions of automorphic representations}
\label{subsec:L-function}

The $L$-function\index{L-function, $L(s, \pi)$} associated to $\pi = \otimes_v \pi_v \in \mathcal{A}(G)$ is of the form
\begin{equation}
\label{L-function}
L(s, \pi) = \prod_{\p} L(s, \pi_\p) = \sum_{N\q \geqslant 1} \frac{a_\pi(\q)}{N\q^s},
\end{equation}

\noindent for $s$ of sufficiently large real part, where $a_\pi(\q)$ are complex numbers, the sum runs over nonzero integer ideals $\q$ of $\mathcal{O}$ and $N\q$ denotes the norm of~$\q$. Here, the $L(s, \pi_\p)$ are the local factors associated to the components $\pi_\p$ at finite places $\p$, and can be written as
\begin{equation}
\label{f-L-factors}
L(s, \pi_\p) = \left(1-\alpha_\pi(\p)N\p^{-s}\right)^{-1}\left(1-\beta_\pi(\p)N\p^{-s}\right)^{-1},
\end{equation}

\noindent where $\alpha_\pi(\p)$ and $\beta_\pi(\p)$ are complex numbers\label{alpha}, called spectral parameters\index{Satake parameters} of~$\pi$ and generalizing the usual Satake parameters for unramified representations. For archimedean places, there are also complex numbers still denoted $\alpha_{\pi}(v)$ and $\beta_{\pi}(v)$ such that the associated local factors take the form
\begin{equation*}
\label{a-L-factors}
L(s, \pi_v) =  \Gamma_v(s-\alpha_{\pi}(v))\Gamma_v(s-\beta_{\pi}(v)),
\end{equation*}

\noindent where $\Gamma_v(s)$ is defined by $\Gamma_\R(s) = \pi^{-s/2}\Gamma(s/2)$ when $v$ is a real place, and by $\Gamma_\C(s) = \Gamma_\R(s)\Gamma_\R(s+1)$ when $v$ is a complex place. The product of these archimedean $L$-factors is denoted $L(s, \pi_\infty)$ and called the archimedean part of the $L$-function. Introduce $c(\pi_f)$ the arithmetic conductor of $\pi \in \mathcal{A}(G)$ and define the completed $L$-function as
\begin{equation*}
\label{completed-L}
\Lambda(s, \pi) = c(\pi_f)^{s/2} L(s, \pi_\infty) L(s, \pi).
\end{equation*}

\noindent It satisfies the functional equation
\begin{equation*}
\label{fe}
\Lambda(s, \pi) = \varepsilon_\pi \Lambda(1-s, {\pi}),
\end{equation*}

\noindent where $\varepsilon_\pi$ is the root number\index{root number, $w_\pi$} of $\pi$ and is among $1$ and $-1$ since $\pi$ is self-dual.


\subsection{Zeros and spectral parameters}
\label{sec:ef}
\label{subsec:ef}


Explicit formulas\index{explicit formula} relate zeros of an $L$-function and prime numbers and thus are a relevant tool to handle one-level densities. The explicit formula of Rudnick and Sarnak \cite{rudnick_zeroes_1996} is in this case particularly well-suited. Their result is stated for the base field $\Q$, but carries on to the setting of more general number fields without particular difficulty, so that we can state a generalized version.

\begin{prop}[Explicit formula]
Let $\pi \in \mathcal{A}(G)$. For all $\phi~\in~\mathcal{S}(\R)$ and $R > 0$,
\begin{align*}
& \sum_{\rho_\pi} \phi\left( \tilde{\gamma}_\pi \right) = \widehat{\phi}(0) \frac{\log c(\pi)}{\log R} \\
& \quad - 2 \sum_\p \sum_{\nu \geqslant 1}(\alpha_\pi^\nu(\p)+\beta^\nu_\pi(\p)) \widehat{\phi}\left(\frac{\nu \log N\p}{\log R} \right) \frac{\log N\p}{N\p^{\nu / 2} \log R} + O\left(\frac{1}{\log R}\right)
\end{align*}

\noindent where the sum on the left hand side runs through the zeros $\rho_\pi = \frac{1}{2} + i\gamma_\pi$ of $L(s, \pi)$.
\end{prop}

\textit{Proof.} The zeros of $\Lambda(s,\pi)$ weighted by their multiplicities correspond to the poles of $\Lambda'(s, \pi)/\Lambda(s, \pi)$ weighted by their residues, and the explicit formula comes from a double evaluation of the integral
\begin{equation}
\label{omg}
I = \frac{1}{2i\pi} \int_{\Re(s) = 2} \frac{\Lambda'}{\Lambda}(s, \pi) {\phi}(s) \dd s.
\end{equation} 

The $L$-function decomposes as an Euler product $\Lambda = \prod_v L_v$, where $L_v(s, \pi)$ stands for $L(s, \pi_v)$ for simplicity. Thus, integrating its logarithmic derivative leads to
\begin{equation}
\label{integral-fe}
I = \sum_v \frac{1}{2i\pi} \int_{\Re(s) = 2} \frac{L_v'}{L_v}(s, \pi) \phi(s) \dd s .
\end{equation}

Denote $I_v$ the integrals appearing in the sum above, and first consider the finite places. Since $L_\p(s, \pi) = (1-\alpha_\pi(\p)N\p^{-s})^{-1} (1-\beta_\pi(\p)N\p^{-s})^{-1}$ at every finite place, the quotient appearing in $I_\p$ can be rewritten
\begin{align*}
\frac{L_\p'}{L_\p}(s, \pi) & = - \frac{\alpha_\pi(\p) N\p^{-s}}{1-\alpha_\pi(\p)N\p^{-s}} \log N\p - \frac{\beta_\pi(\p) N\p^{-s}}{1-\beta_\pi(\p)N\p^{-s}} \log N\p \\
& = - \log(N\p) \sum_{\nu \geqslant 1} (\alpha_\pi(\p)^\nu + \beta_\pi(\p)^\nu) N\p^{-\nu s}.
\end{align*}

Let $\phi^\star(x) = \phi(1/2 + x)$. Since $\phi$ is an holomorphic function, it has no poles and thus the contour appearing in the integral \eqref{integral-fe} can be translated from the vertical line of abscissa $2$ to the one of abscissa $\frac{1}{2}$. Introduce the inverse Mellin transform $\mathcal{M}\phi$ of $\phi$. Using change of variables, we get
\begin{align*}
\mathcal{M} \phi (y) & = \frac{1}{2i\pi}\int_{\Re (s) = 1/2} \phi(s) y^{-s} \dd s \\
& =  \frac{y^{-1/2}}{2\pi} \int_\R \phi\left(\frac{1}{2} + ir\right) e^{-i r \log y} \dd r \\
& = \frac{y^{-1/2}}{2\pi} \widehat{\phi^\star} (\log y),
\end{align*}

\noindent so that the finite local integrals become
\begin{align*}
I_\p & = - \log(N\p) \sum_{\nu \geqslant 1} (\alpha_\pi(\p)^\nu + \beta_\pi(\p)^\nu) \mathcal{M} \phi (N \p^\nu) \\
& = - \frac{1}{2\pi}\sum_{\nu \geqslant 1} (\alpha_\pi(\p)^\nu + \beta_\pi(\p)^\nu) \widehat{\phi^\star} \left( \nu \log N\p \right) \frac{\log N\p}{N\p^{\nu/2}}.
\end{align*}

\noindent  If $\phi$ is $\frac{1}{2}$-symmetrical, then letting $\phi^\vee (s) = \phi(1-s)$ and applying the computations above, we have
\begin{equation*}
\label{Mellin-symmetry}
\frac{y^{-1/2}}{2\pi} \widehat{\phi^\star} (- \log y) = \frac{y^{-1/2}}{2\pi} \int_\R \phi\left(\frac{1}{2} - ir\right) e^{-i r \log y} \dd r = \mathcal{M} \phi^\vee (y) =  \mathcal{M} \phi (y).
\end{equation*}

On the other hand, the Cauchy theorem allows to unfold the integral \eqref{omg} in terms of the zeros of $\Lambda(s, \pi)$. Indeed, since ${\phi}$ is entire, the only poles of the integrated function are the zeros of $\Lambda$ and the corresponding residues are their multiplicities. All these zeros lie in the vertical strip $-1 < \Re(s) < 2$, so that translating the contour through this whole band captures all the zeros and gives
\begin{equation}
I = \sum_{\rho_\pi} \phi(\rho_\pi) + \frac{1}{2i\pi} \int_{\Re(s) = -1} \frac{\Lambda'}{\Lambda}(s, \pi) \phi(s) \dd s .
\end{equation}

The functional equation of $L$ is of the form $\Lambda(s, \pi) = \varepsilon_\pi \Lambda(1-s, \pi)$, thus inputing it in the integral above and changing variables to come back to the vertical line of abscissa $\Re(s) = 2$, we get
\begin{equation}
I = \sum_{\rho_\pi} \phi(\rho_\pi)  -  \frac{1}{2i\pi} \int_{\Re(s) = 2} \frac{\Lambda'}{\Lambda}(s, \pi) \phi(1-s) \dd s .
\end{equation}

Coming back to the definition \eqref{integral-fe} of $I$,  we deduce
\begin{equation}
\sum_{\rho_\pi} \phi(\rho_\pi) = \frac{1}{2i\pi} \int_{\Re(s) = 2} \frac{\Lambda'}{\Lambda}(s, \pi) (\phi(1-s)+\phi(s)) \dd s .
\end{equation}

Finally, for an even Schwartz function $\phi_0$, the function 
\begin{equation}
\phi(s) = \phi_0 \left( \frac{\log R}{2\pi} \left( s - \frac{1}{2}\right) \right), \qquad R > 0,
\end{equation}

\noindent is Schwartz and $\frac{1}{2}$-symmetric, and moreover satisfies
\begin{equation}
\widehat{\phi^\star}(\nu \log N\p) = \frac{2\pi}{\log R} \widehat{\phi}_0 \left( \frac{\nu \log N\p}{\log R} \right).
\end{equation}

The archimedean places are dealt with by the same the treatment as in Rudnick and Sarnak with no modification, since the associated archimedean local fields remain among $\R$ and $\C$. It contributes as an error term plus a factor carrying the archimedean part of the conductor. Combining the two previous expressions of the integral,
\begin{align*}
& \sum_{\rho_\pi} \phi(\rho_\pi)  = \frac{\log c(\pi) }{2\pi} \int_{\Re(s) = 2} \phi(s) \dd s\\
& \quad - \frac{1}{\pi} \sum_\p  \sum_{\nu \geqslant 1} (\alpha_\pi(\p)^\nu + \beta_\pi(\p)^\nu) \widehat{\phi^\star} \left( \nu \log N\p \right) \frac{\log N\p}{N\p^{\nu/2}} + O\left(\frac{1}{\log R}\right) \\
&   = \frac{\log c(\pi) }{\log R} \widehat{\phi}_0(0) \\
& \quad - 2 \sum_\p  \sum_{\nu \geqslant 1} (\alpha_\pi(\p)^\nu + \beta_\pi(\p)^\nu) \widehat{\phi}_0 \left( \frac{\nu \log N\p}{\log R} \right) \frac{\log N\p}{N\p^{\nu/2} \log R} + O\left(\frac{1}{\log R}\right)
\end{align*}

\noindent ending the generalization of the result of Rudnick and Sarnak. \qed

\noindent \textit{Remark.} Two different objects are deliberately denoted by $\pi$ in the proof above, namely an automorphic representation and the pythagorean constant. Convinced that no confusion should arise, we preferred to stick with the long lasting tradition of both notations.

The explicit formula then yields, with $R = c(\pi)$, the reformulation of the one-level density
\begin{align}
\label{ef}
& D(\pi, \phi)= \widehat{\phi}(0) \\
& \quad - \frac{2}{\log c(\pi)} \sum_{\p} \sum_{\nu = 1}^\infty \left(\alpha_\pi^\nu(\p) + \beta_\pi^\nu(\p)\right) \widehat{\phi}\left(\frac{\nu\log N\p}{\log c(\pi)}\right)\frac{\log N\p}{N\p^{\nu/2}} + O\left(\frac{1}{\log c(\pi)}\right).\notag
\end{align}

After switching summations, consider the inner sum for a fixed $\nu \geqslant 1$, that is to say
\begin{equation}
\label{Pnu-pi}
P^{(\nu)}(\pi, \phi) = \frac{2}{\log c(\pi)} \sum_{\p} (\alpha_\pi^\nu (\p)+ \beta_\pi^\nu(\p)) \widehat{\phi}\left(\frac{\nu \log N\p}{\log c(\pi)}\right)\frac{\log N\p}{{N\p}^{\nu / 2}},
\end{equation}

\noindent so that the one-level density decomposes as
\begin{equation}
\label{old}
D(\pi, \phi) = \widehat{\phi}(0) - \sum_{\nu \geqslant 1} P^{(\nu)}(\pi , \phi) + O\left(\frac{1}{\log c(\pi)}\right).
\end{equation}

For a given automorphic representation $\pi$, there is only a finite number of zeros of $L(s, \pi)$ in the compact support of $\widehat{\phi}$. Following the philosophy of Katz and Sarnak, introduce rather the average over the truncated universal family
\begin{equation}
\label{P-averaged}
{\mathcal{P}}^{(\nu)}_{Q}(\phi)  = \frac{1}{|\mathcal{A}(Q)|} \sum_{\pi \in \mathcal{A}(Q)} P^{(\nu)}(\pi, \phi).
\end{equation}

\noindent The following sections are dedicated to estimate the contribution of these ${\mathcal{P}}^{(\nu)}_{Q}$.

\section{Decomposition of the universal family}
\label{sec:decomposition}

\subsection{Sieving the universal family}

In order to prove Theorem \ref{thmD}, it is necessary to decompose the universal family into smaller sets with fixed spectral data, amenable to trace formula methods. The conductor of $\pi \in \mathcal{A}(G)$ splits into local conductors by its very definition, in particular can be written
\begin{equation}
\label{conductor-splitting}
c(\pi) = c(\pi_R)   c(\pi^R).
\end{equation}

\noindent Thus, the universal family decomposes as
\begin{equation}
\label{uf-decomposition}
\mathcal{A}(Q) =  \bigsqcup_{\substack{N\q \leqslant Q\\ \q \wedge R = 1}} \bigsqcup_{\substack{\sigma_R \in \widehat{G}_R \\ c(\sigma_R) \leqslant Q/N\q}}  \mathcal{A}(\q, \sigma_R)  ,
\end{equation}

\noindent where the sets $\mathcal{A}(\q, \sigma_R) $ are obtained by fixing spectral data to an arithmetic conductor $\q$ at split places and to a ramified isomorphism class $\sigma_R \in \widehat{G}_R$, that is
\begin{align*}
\mathcal{A}(\q, \sigma_R) & =  \left\{ \pi \in \mathcal{A}(G) \ : \ \pi_R \simeq \sigma_R, \ \c(\pi^{R}) = \q \right\}.
\end{align*}

\noindent This decomposition of the universal family reduces the study of $\mathcal{A}(G)$ to the harmonic subfamilies $\mathcal{A}(\q, \sigma_R)$, easier to grasp in the context of trace formulas. The crucial point is to replace the global condition of belonging to $\mathcal{A}(Q)$ by local conditions. For $\phi$ a Schwartz class function on $\R$ with compactly supported Fourier transform, the above partition induces a decomposition of the ${\mathcal{P}}_Q^{(\nu)}$ as
\begin{equation}
\label{decomposition1}
\begin{split}
{\mathcal{P}}_Q^{(\nu)}(\widehat{\phi}\,) & =  \frac{1}{Q^2} \sum_{\pi \in \mathcal{A}(Q)} P^{(\nu)}(\pi, \phi)  \\
& =  \frac{1}{Q^2} \sum_{\substack{\pi \in \mathcal{A}(G) \\ c(\pi_R)c(\pi^R) \leqslant Q}} P^{(\nu)}(\pi, \phi)  \\
& = \frac{1}{Q^2}  \sum_{\substack{\sigma_R \in \widehat{G}_R \\ c({\sigma_R}) \leqslant Q}}  \sum_{\substack{N\q \leqslant Q/c(\sigma_R) \\ \q \wedge  R = 1}} 
\sum_{\substack{\pi \in \mathcal{A}(\q, \sigma_R) }} P^{(\nu)}(\pi, \phi)  
\end{split}
\end{equation}

\noindent where the sum over $\q$ is meant to run through ideals of $\mathcal{O}^{R}$. Introduce~${\mathcal{P}}^{(\nu)}_{\q, \sigma_R}(\phi)$ the innermost parts of the splitting in the first summation above, that is
\begin{equation}
\label{A}
{\mathcal{P}}^{(\nu)}_{\q, \sigma_R}(\PPP) = \sum_{\substack{\pi \in \mathcal{A}(\q, \sigma_R)}} P^{(\nu)}(\pi, \phi) .
\end{equation}

\noindent \textit{Remark.} In the case of non totally definite quaternion algebras, this decomposition needs to be refined by using the archimedean Langlands classification of the unitary spectrum in order to take into account the continuous spectrum at archimedean split places, see \cite{brumley_counting_2016} and \cite{lesesvre_counting_2020}.

Expanding the expression of $P^{(\nu)}(\pi, \phi)$ given by the explicit formula and switching summations lead to
\begin{equation}
\label{Pnuq}
{\mathcal{P}}^{(\nu)}_{\q, \sigma_R}(\phi) = \sum_\p \left( \sum_{\pi \in \mathcal{A}(\q, \sigma_R)} \left( \alpha_\pi^\nu(\p) + \beta_\pi^\nu(\p)  \right) \right) \widehat{\phi}\left( \frac{\nu \log N\p}{\log c(\pi)} \right) \frac{2\log N\p}{N\p^{\nu / 2} \log c(\pi)},
\end{equation}
\noindent where $c(\pi)$ stands as a shortcut notation for $N\q c(\sigma_R)$, justifying its presence outside the sum over the harmonic subfamily $\mathcal{A}(\q, \sigma_R)$ in which the conductor is restricted to $\c(\pi^R) = \q$ and $\pi_R \simeq \sigma_R$. This convention will be steadily used in the following. Introduce the spectral sums,
\begin{align}
\label{lambda}
\tilde{\Lambda}^{(\nu)}_{\q, \sigma_R}(\p)  = \sum_{\pi \in \mathcal{A}(\q, \sigma_R)} \left( \alpha_\pi^\nu(\p) + \beta_\pi^\nu(\p)  \right).
\end{align}

\subsection{Old and new forms}

The one-level density \eqref{density-llz} sees no multiplicities\index{multiplicities}, but the trace formula counts them. The spectral multiplicities associated to the decomposition of $L^2(G(F) \backslash G(\A))$, which are more suitable weights for the forthcoming computations, are given by 
\begin{equation}
\label{multiplicities}
m\left(\pi, \q\right) = \dim \ \pi^{\overline{K}_0(\q)},
\end{equation}

\noindent where
\begin{equation}
Z K_0(\q) = \prod_{\substack{\p^r |\!| \q}} Z_\p K_{0, \p}\left(\p^r\right) \subseteq B^{\times}\left(\A^R_f\right),
\end{equation}

\noindent and $\overline{K}_0(\q)$ stands for the image of $ZK_0(\q)$ under the natural projection $B^\times \rightarrow G$. The choice is made so that $m(\pi, \q) \neq 0$ is equivalent to $\c(\pi^R_f) \div \q$. The analogous sum to \eqref{lambda} weighted by the multiplicities is
\begin{equation}
\label{B}
{\Lambda}^{(\nu)}_{\q, \sigma_R}(\p) = \sum_{\substack{\pi \in \mathcal{B}(\q, \sigma_R)}} m\left(\pi, \q\right) \left( \alpha_\pi^\nu(\p) + \beta_\pi^\nu(\p) \right), 
\end{equation}

\noindent where
\begin{equation}
\mathcal{B}(\q, \sigma_R) = \left\{ \pi \in \mathcal{A}(G) \ : \ \pi_R \simeq \sigma_R, \  \c\left(\pi^{R}_f\right) \div \q \right\}.
\end{equation}

The sum defined by \eqref{A} runs over the newforms\index{newforms} while \eqref{B} runs over the old ones. The relation between them lies in the following lemma.

\begin{lemma}  
\label{lem:sieve}
Let $\q$ be an integer ideal prime to $R$ and $\sigma_R$ an irreducible unitary representation of $G_R$. Let $\lambda_2= \mu \star \mu$ where $\mu$ is the Möbius function. For every $Q \geqslant 1$,
\begin{align*}
\tilde{\Lambda}^{(\nu)}_{\q, \sigma_R}(\p)  & = \sum_{\d\div \q} \lambda_2\left(\frac{\q}{\d}\right) {\Lambda}^{(\nu)}_{\d, \sigma_R}(\p)  .
\label{sieve}
\end{align*}
\end{lemma}

\proof Recall that, for every finite split place $\p$, Casselman gives the local multiplicites of $\sigma_\p \in \widehat{G}_\p$,
\begin{equation}
\dim \sigma_{\p}^{K_0\left(\p^{\f(\sigma_\p)+i}\right)} = i+1, \qquad i \geqslant 0.
\end{equation}

\noindent From this immediately follows for a $\sigma \in \mathcal{A}(G)$, after taking the product over all finite split places, that the global multiplicities are
\begin{equation}
m\left(\sigma, \q\right) = \tau_2\left(\frac{\q}{\c(\sigma^{R})}\right),
\label{mult}
\end{equation}

\noindent where $\tau_2 = 1 \star 1$ is the divisor function. Since $\left(\sigma^{R}\right)^{\overline{K}_0\left(\q\right)} \neq 0$ implies $\c(\sigma^{R}) \div \q$, the sum defining ${\Lambda}^{(\nu)}_{\q, \sigma_R}(\p) $ is eventually reduced to a sum over $\c(\sigma^{R}) \div \q$. Thus, by the precise knowledge \eqref{mult} of the multiplicities, 
\begin{equation}
\begin{split}
{\Lambda}^{(\nu)}_{\q, \sigma_R}(\p)  &  = \sum_{\d\div \q} \sum_{\sigma \in \mathcal{A}(\d, \sigma_R)} \tau_2\left(\frac{\q}{\c(\sigma^{R})}\right) P^{(\nu)}(\sigma, \phi)  \\
&  = \sum_{\d\div \q} \tau_2\left(\frac{\q}{\d}\right) \sum_{\sigma \in \mathcal{A}(\d, \sigma_R)} P^{(\nu)}(\sigma, \phi) \\
&  = \sum_{\d\div \q} \tau_2\left(\frac{\q}{\d}\right) \tilde{\Lambda}^{(\nu)}_{\d, \sigma_R}(\p) 
\end{split}
\end{equation}

\noindent so that $\Lambda = \tau_2 \star \tilde{\Lambda}$, with a slight abuse of notation. Hence, by Möbius inversion,
\begin{equation}
\tilde{\Lambda}^{(\nu)}_{\q, \sigma_R}(\p)  = \sum_{\d\div \q} \lambda_2\left(\frac{\q}{\d}\right) {\Lambda}^{(\nu)}_{\d, \sigma_R}(\p) ,
\end{equation}

\noindent achieving the proof of the claim. \qed

So that according to the decomposition of the universal family~\eqref{uf-decomposition}, the average \eqref{P-averaged} can be rewritten
\begin{equation}
\label{sieve}
{\mathcal{P}}^{(\nu)}_{Q}(\phi) = \frac{1}{|\mathcal{A}(Q)|}\sum_{\substack{N\q \leqslant Q \\ \q \wedge R = 1}} \sum_{\substack{\sigma_R \in \widehat{G}_R \\ c(\sigma_R) \leqslant Q/N\q}} \mathcal{P}^{(\nu)}_{\q, \sigma_R}(\phi),
\end{equation}

\noindent where we rewrite
\begin{equation*}
{\mathcal{P}}^{(\nu)}_{\q, \sigma_R}(\phi) = \sum_\p \left( \sum_{\d | \q} \lambda_2\left(\frac{\q}{\d}\right) \Lambda_{\d, \sigma_R}^{(\nu)}(\p) \right) \widehat{\phi}\left( \frac{\nu \log N\p}{\log c(\pi)} \right) \frac{2\log N\p}{N\p^{\nu / 2} \log c(\pi)}.
\end{equation*}

\section{High orders contributions}
\label{subsec:high-contrib}
\label{sec:high-contrib}

For $\nu$ large enough, it is possible to bound directly $P^{(\nu)}(\pi, \phi)$ and show that they do not contribute to the type of symmetry.

\begin{prop}
\label{prop:high-contrib}
For $Q \geqslant 1$,
\begin{equation}
\label{high-contrib}
\sum_{\nu \geqslant 3} {\mathcal{P}}^{(\nu)}_{Q}(\phi) \ll \frac{1}{\log Q}.
\end{equation}
\end{prop}
	
\proof The main aim is to bound the spectral parameters $\alpha_\pi^\nu(\p) + \beta_\pi^\nu(\p)$ in the sum \eqref{Pnu-pi}. For holomorphic cusp forms, the Ramanujan conjecture\index{Ramanujan conjecture} holds by Deligne \cite{deligne_conjecture_1974} and states that $|\alpha_\pi(\p) + \beta_\pi(\p)| \leqslant 2$. For Maass forms over general number fields, Blomer and Brumley \cite{blomer_ramanujan_2011} proved that 
\begin{equation}
|\alpha_\pi(\p) + \beta_\pi(\p)| \ll N\p^{7/64}.
\end{equation}

\noindent Hence for any cuspidal automorphic representation of $\GL(2)$ this last bound is valid, in particular for $\pi \in \mathcal{A}(G)$. Thus,
\begin{align*}
\sum_{\nu \geqslant 3} P^{(\nu)}(\pi, \phi) & \ll \frac{1}{\log c(\pi)} \sum_{\p} \sum_{\nu \geqslant 3} \frac{\log N\p}{N\p^{\nu(1/2 - 7/64)}}\\
& \ll \frac{1}{\log c(\pi)} \sum_{\p} \frac{\log N\p}{N\p^{3\left(1/2-7/64 \right)}} \\
& \ll \frac{1}{\log c(\pi)}
\end{align*}

\noindent Turning back to the sums over partial families, introducing the cardinality $A(\q, \sigma_R)$ of $\mathcal{A}(\q, \sigma_R)$ and using the fact that $c(\pi) \geqslant N\q$,  we get
\begin{equation}
\label{AZE}
\sum_{\nu \geqslant 3} {\mathcal{P}}^{(\nu)}_{\q, \sigma_R}(\phi) \ll  \sum_{\pi \in \mathcal{A}(\q, \sigma_R)} \frac{1}{\log N\q}  \ll \frac{A(\q, \sigma_R)}{\log N\q}.
\end{equation}

The cardinality of the family $\mathcal{A}(\q, \sigma_R)$ has been computed in our previous work \cite{lesesvre_counting_2020}: there is a remainder term $R(\q, \sigma_R)$ such that
\begin{equation}
\label{cardinalityP}
A(\q, \sigma_R) = \vol(G(F)\backslash G(\A))\varphi_2(\q)\mu^\Pl(\sigma_R) + R(\q, \sigma_R),
\end{equation}

\noindent where  $\varphi_2 =  \lambda_2 \star \mu^2 \star \id$ and the remainder satisfies, for a certain $\theta > 0$, 
\begin{equation}
\label{remainder-R}
\sum_{\substack{\q \leqslant Q \\ \q \wedge R = 1}} \sum_{\substack{\sigma_R \in \widehat{G}_R \\ c(\sigma_R) \leqslant Q/N\q}} R(\q, \sigma_R) \ll Q^{2-\theta}.
\end{equation}


Introduce the following dampening lemma, justifying that the logarithmic factor appearing in the denominator of \eqref{AZE} is enough to make the whole sum negligible compared to the one free of this factor, that is to say the cardinality of the truncated universal family $|\mathcal{A}(Q)|$.

\begin{lemma}
\label{lem:dampening}
Let $f$ be a positive function on the integer ideals of $F$, for which there is an $\alpha >0$ such that
\begin{equation}
\sum_{N\n \leqslant X} f(\n) \sim X^\alpha.
\end{equation}

\noindent Then for every $\varepsilon > 0$,
\begin{equation}
\label{assumption}
\sum_{N\n \leqslant X} \frac{f(\n)}{\log(N\n)^\varepsilon} \ll  \frac{1}{\log(X)^\varepsilon} \sum_{N\n \leqslant X} f(\n).
\end{equation}
\end{lemma}

\proof The hyperbola method can be efficiently used in this setting. Cutting the sum at $X^{1/2}$ for a positive $X$ leads to
\begin{align*}
\sum_{N\n \leqslant X} \frac{f(\n)}{\log(N\n)^\varepsilon}  & = \sum_{N\n \leqslant X^{1/2}} \frac{f(\n)}{\log(N\n)^\varepsilon} + \sum_{X^{1/2} < N\n \leqslant X} \frac{f(\n)}{\log(N\n)^\varepsilon} \\
& \ll \sum_{N\n \leqslant X^{1/2}} f(\n) + \frac{1}{\log(X)^{\varepsilon}} \sum_{X^{1/2} < N\n \leqslant X} f(\n) \\
& \ll \frac{1}{\log(X)^\varepsilon}  \sum_{N\n \leqslant X} f(\n).
\end{align*}

\noindent Indeed, the asymptotic assumption \eqref{assumption} yields
\begin{equation}
\sum_{N\n \leqslant X^{1/2}} f(\n) \ll X^{\alpha /2} \ll  \frac{X^{\alpha}}{\log(X)^\varepsilon} \ll \frac{1}{\log(X)^\varepsilon}  \sum_{N\n \leqslant X} f(\n),
\end{equation}

\noindent proving the lemma. \qed

So after summations of the contributions of the sieved families \eqref{cardinalityP}, Lemma \ref{lem:dampening} provides the bound
\begin{align*}
 \sum_{\nu \geqslant 3} {\mathcal{P}}_Q^{(\nu)}(\phi) & = \frac{1}{Q^2}\sum_{\substack{N\q \leqslant Q \\ \q \wedge R = 1}} \sum_{\substack{\sigma_R \in \widehat{G}_R \\ c(\sigma_R) \leqslant Q/N\q}} \sum_{\nu \geqslant 3} {\mathcal{P}}_{\q, \sigma_R}^{(\nu)}(\phi) \\
& \ll  \frac{1}{Q^2}\sum_{\substack{\sigma_R \in \widehat{G}_R \\ c(\sigma_R) \leqslant Q}} \sum_{\substack{N\q \leqslant Q/c(\sigma_R) \\ \q \wedge R = 1}}  \frac{A(\q, \sigma_R)}{\log N\q} \\
& \ll \frac{1}{Q^2} \sum_{\substack{\sigma_R \in \widehat{G}_R \\ c(\sigma_R) \leqslant Q}} \mu^{\mathrm{Pl}}(\sigma_R) \sum_{\substack{N\q \leqslant Q/c(\sigma_R) \\ \q \wedge R = 1}} \frac{\varphi_2(\q)}{\log N\q} \\
& \qquad +\frac{1}{Q^2} \sum_{\substack{\sigma_R \in \widehat{G}_R \\ c(\sigma_R) \leqslant Q}} \sum_{\substack{N\q \leqslant Q/c(\sigma_R) \\ \q \wedge R = 1}}  R(\q, \sigma_R) \\
& \ll \frac{1}{\log Q}\sum_{\substack{\sigma_R \in \widehat{G}_R \\ c(\sigma_R) \leqslant Q}} \frac{\log c(\sigma_R)}{c(\sigma_R)^2} \mu^{\mathrm{Pl}}(\sigma_R) + Q^{-\theta}
\end{align*}

\noindent and the sum over the ramified spectrum converges by \cite{lesesvre_counting_2020}. This ends the proof of Proposition \ref{prop:high-contrib}. \qed

\rk Analogously to the work of Iwaniec, Luo and Sarnak \cite{iwaniec_low_2000} and most of the literature on low-lying zeros, the high order terms are negligible, with a logarithmic savings. This bound follows from directly dominating $P^{(\nu)}(\pi, \phi)$ without use of neither the average over the family nor the sum over the primes.

\section{Hecke operators}
\label{sec:traces-hecke}

\subsection{Hecke eigenvalues and coefficients}
\label{subsec:hecke-opertaors-def}
\label{subsec:hecke-opertaors-eigenvalues}
\label{subsec:hecke-operators-relation-coeff-eigenvalues}


For the two remaining cases $\nu = 1$ and $\nu=2$, straightforward estimations are no more sufficient, feature already present in the axiomatic proposed by Dueñez and Miller \cite{duenez_effect_2009} for the behavior of low-lying zeros in families of L-functions. The inner spectral sums
\begin{equation}
\sum_{\pi \in \mathcal{B}(\q, \sigma_R)} m(\pi, \q) \left( \alpha_\pi^\nu(\p) + \beta_\pi^\nu(\p) \right)
\end{equation}

\noindent are closely related to traces of Hecke operators, so that $\mathcal{P}^{(\nu)}_{\q, \sigma_R}(\phi)$ should be interpreted as a spectral side of a trace formula.  Define the normalized Hecke operators \index{Hecke operator, $T_{\p^\nu}$} as
\begin{equation}
\label{def:hecke}
T_{\p^\nu} = N\p^{-\nu/2} \mathbf{1}_{T(\p^\nu)} \quad \text{where} \quad T(\p^\nu) = 
\bigcup_{\substack{i+j = \nu \\ 0 \leqslant i \leqslant j}} 
K_\p
\left(
\begin{array}{cc}
\p^i & \\ & \p^j
\end{array}
\right)
K_\p.
\end{equation}

The Hecke operator for a global ideal $\mathfrak{n}$ of $\mathcal{O}$ is defined by
\begin{equation}
T_{\n} = \prod_{\p^r |\!| \mathfrak{n}} T_{\p^r}.
\end{equation}

One of the main appeal of Hecke operators is that they provide an explicit recipe to catch the coefficients of $L$-functions. Indeed, they satisfy the same induction relation, hence are equal once suitably normalized.  This is the content of the following standard proposition. Recall that the coefficients of the Dirichlet series attached to an automorphic representation are normalized no that $a_\pi(1) = 1$.

\begin{prop}
\label{prop:hecke-eigenvalues-equal-coefficients}
Let $\mathfrak{n}$ be an ideal of $\mathcal{O}$ and $\pi$ be an unramified representation at the places dividing $\mathfrak{n}$. Introduce $\lambda_\pi(\n)$ the eigenvalue of~$T_{\n}$ acting on $\pi$. Then , 
\begin{equation}
a_\pi(\n) = \lambda_\pi(\n).
\end{equation}

\noindent Moreover, if $\pi$ is ramified at one of the places dividing $\mathfrak{n}$, then $T_{\mathfrak{n}}$ acts by zero on $\pi$.
\end{prop}

\proof The $T_{\p^n}$ satisfy \cite[Prop 4.6.4]{bump_automorphic_1997} the recursive relation
\begin{equation}
\label{Hecke-rec-relation}
T_{\p^{n+1}} = T_\p T_{\p^n} -  T_{\p^{n-1}}, \qquad n \geqslant 1,
\end{equation}

\noindent which transfers at the Hecke eigenvalues\index{Hecke eigenvalue, $\lambda_\pi(\p^n)$} level and gives
\begin{equation}
\label{Hecke-ev-rec-relation}
\lambda_\pi(\p^{n+1}) = \lambda_\pi(\p) \lambda_\pi(\p^n) - \lambda_\pi(\p^{n-1}), \qquad n \geqslant 1.
\end{equation}

Recall that, by the Euler product decomposition of $L(s, \pi)$, the coefficients $a_\pi(\n)$ are entirely determined by their values at the prime powers~$\p^n$. Moreover, the centerless setting implies a trivial central character and hence the Satake parameters are related by $\alpha_\pi(\p) = \beta_\pi(\p)^{-1}$. Let $\alpha$ be a shortened version of $\alpha_\pi(\p)$. The local $L$-factors are
\begin{equation*}
\label{Hecke-L-dvt}
L_\p(s,\pi) = (1-\alpha N\p^{-s})^{-1} (1-\alpha^{-1} N\p^{-s})^{-1}
= \left(\sum_{i \geqslant 0} \alpha^i N\p^{-is} \right)\left( \sum_{j \geqslant 0} \alpha^{-j} N\p^{-js} \right).
\end{equation*}

By unfolding the power series, the coefficient of $N\p^{ns}$ is 
\begin{equation}
\label{coefficients-L-series}
a_\pi(\p^n) = \sum_{i+j=n} \frac{\alpha^i}{\alpha^j}.
\end{equation}

A straightforward computation hence leads to the recursion relation
\begin{equation}
\label{coefficients-relation}
a(\p^{n+1}) = a(\p)a(\p^n) - a(\p^{n-1}).
\end{equation}

Since this is the same relation than for the Hecke eigenvalues \eqref{Hecke-rec-relation}, the two sequences are proportional. Moreover, $\p$-unramified newforms are normalized so that $a_\pi(1) = \lambda_\pi(1)=1$, leading to the equality of both sequences as claimed. The second part of the proposition is straightforward, since the~$T_\p$ are the unramified Hecke operators, in particular are bi-$K_\p$-invariants, and so project on $K_\p$ fixed vectors, which are reduced to zero for $\p$-ramified representations. \qed

\subsection{Spectral selection}

The above proposition states that Hecke eigenvalues are a way to interpret the coefficients of an automorphic representation. Since these coefficients are related to the Satake parameters by the Euler product development, this fact provides a way to handle the sums $\Lambda_{\q, \sigma_R}^{(\nu)}(\p)$. In order to make this question amenable to the trace formula method, it is necessary to have a grasp on these quantities by means of Fourier transforms. This is provided by the following standard proposition.

\begin{lemma}
\label{lem:hecke}
\label{prop:hecke-eigenvalues}
Every $K_\p$-spherical vector $v$ in a representation $\pi_\p$ is an eigenvector for all the Hecke operators $T_\p^n$. Moreover, for every $\nu \geqslant 0$,
\begin{equation}
\label{hecke-action}
\widehat{T}_{\p^\nu} (\pi_\p) =
\left\{
\begin{array}{cl}
\lambda_\pi(\p^\nu) & \text{if } \pi_\p \text{ is unramified at } \p; \\
0 & \text{otherwise.}
\end{array}
\right.
\end{equation}
\end{lemma}


The appearing dissymmetry  between the ramified and unramified cases motivates the introduction of
\begin{align}
\label{specsum, Lambda}
\Lambda^{(\nu), \ur}_{\q, \sigma_R}(\p) & = \sum_{\pi \in \mathcal{B}(\q, \sigma_R)} m(\pi, \q) \lambda_\pi(\p^\nu)\\
\label{specsumr}
\Lambda^{(\nu), \r}_{\q, \sigma_R}(\p) & = \sum_{\pi \in \mathcal{B}(\q, \sigma_R)^{\mathrm{r}, \p}} m(\pi, \q) \left( \alpha_\pi^\nu(\p) + \beta_\pi^\nu(\p) \right)
\end{align}

\noindent where $\mathcal{B}(\q, \sigma_R)^{\mathrm{r}, \p}$ denotes the subset of $\mathcal{B}(\q, \sigma_R)$ composed of its representation  ramified at $\p$. Note that the sum in $\Lambda^{(\nu), \ur}_{\q, \sigma_R}(\p)$ is in fact restricted to representations unramified at $\p$ since, by Lemma \ref{lem:hecke}, we have $\lambda_\pi(\p^\nu)=0$ for representations ramified at $\p$. Introduce the analogous notations for $\tilde{\Lambda}^{(\nu), \ur}_{\q, \sigma_R}(\p)$ and $\tilde{\Lambda}^{(\nu), \r}_{\q, \sigma_R}(\p)$.

\subsection{Sums of Hecke eigenvalues}

The previous section showed that Hecke operators lead to a reformulation of the spectral parameters in terms of Hecke eigenvalues in the case of unramified representations. Introduce $\mathcal{P}^{(\nu), \ur}_{\q, \sigma_R}$ to be the contribution of the unramified part of the spectrum to $\mathcal{P}^{(\nu)}_{\q, \sigma_R}(\phi)$, namely
\begin{equation}
\label{blo}
\mathcal{P}^{(\nu), \ur}_{\q, \sigma_R}(\phi) = \sum_{\p} \left( \sum_{\substack{\pi \in \mathcal{B}(\q, \sigma_R) \\ \pi \text{ unram at } \p}} \left( \alpha_\pi^\nu(\p) + \beta_\pi^\nu(\p) \right) \right) \widehat{\phi}\left(\frac{\nu \log N\p}{\log c(\pi)}\right)\frac{2 \log N\p}{{N\p}^{\nu / 2}\log c(\pi)}.
\end{equation}

\subsubsection{Relation in the case $\nu = 1$}

Since $\alpha_\pi(\p) + \beta_\pi(\p) = \lambda_\pi(\p)$ for unramified representations, the corresponding total sum in $\mathcal{P}^{(1), \ur}_{\q, \sigma_R}$ can be rewritten as
\begin{equation}
\label{1-part}
\mathcal{P}^{(1), \ur}_{\q, \sigma_R}(\phi) = \sum_{\p} \tilde{\Lambda}^{(1), \ur}_{\q, \sigma_R}(\p) \widehat{\phi}\left(\frac{\log N\p}{\log c(\pi)}\right)\frac{2 \log N\p}{\sqrt{N\p}\log c(\pi)}.
\end{equation}

\subsubsection{Relation in the case $\nu = 2$}
\label{subsubsec:case}

By identification of the corresponding expressions of the local $L$-factor \eqref{f-L-factors} follows the relation $\alpha^2_\pi(\p) + \beta^2_\pi(\p) = \lambda_\pi(\p^2) - 1$ holding for unramified representations. So that the sum over primes splits as
\begin{align}
\label{P2-split}
\mathcal{P}^{(2), \ur}_{\q, \sigma_R}(\phi) & = \sum_{\p} \tilde{\Lambda}^{(2), \ur}_{\q, \sigma_R}(\p) \widehat{\phi}\left(\frac{2 \log N\p}{\log c(\pi)}\right)\frac{2 \log N\p}{N\p \log c(\pi)} \\
& \quad - \sum_{\p} \widehat{\phi}\left(\frac{2\log N\p}{\log c(\pi)}\right)\frac{2 \log N\p}{N\p\log c(\pi)} .\notag
\end{align}

\noindent The fact that $\widehat{\phi}$ is even and compactly supported in $(-T_\phi, T_\phi)$ allows to write
\begin{equation*}
\int_0^{T_\phi} \widehat{\phi} = \frac{1}{2} \int_\R \widehat{\phi} = \frac{1}{2} \phi(0).
\end{equation*} 

By the Mertens estimates and integration by parts, the second sum in the right hand side rewrites
\begin{align*}
& \sum_{\p} \widehat{\phi}\left(\frac{2\log N\p}{\log c(\pi)}\right)\frac{2 \log N\p}{N\p\log c(\pi)} = \frac{1}{\log c(\pi)} \widehat{\phi} \left( \frac{2\log c(\pi)^{T_\phi / 2}}{\log c(\pi)} \right) \\
& \qquad - \frac{2}{\log c(\pi)} \int_1^{c(\pi)^{T_\phi / 2}} \left( \sum_{\p \leqslant t} \frac{\log N\p}{N\p} \right) \partial_t  \widehat{\phi}\left(\frac{2 \log t}{\log c(\pi)}\right) \dd t \\
& = - \frac{2}{\log c(\pi)} \int_1^{c(\pi)^{T_\phi / 2}} (\log(t)+O(1)) \partial_t \widehat{\phi}\left(\frac{2 \log t}{\log c(\pi)}\right) \dd t\\
& = - \frac{2}{\log c(\pi)} \int_1^{c(\pi)^{T_\phi / 2}} \widehat{\phi}\left(\frac{2 \log t}{\log c(\pi)}\right)  \frac{\dd t}{t} + O\left(\frac{1}{\log c(\pi)}\right) \\
& = \frac{1}{2} \phi(0) + O\left(\frac{1}{\log c(\pi)}\right).
\end{align*}

\noindent The expression \eqref{P2-split} of $\mathcal{P}^{(2), \ur}_{\q, \sigma_R}$ is hence reduced to
\begin{align}
\label{2-part}
\mathcal{P}^{(2), \ur}_{\q, \sigma_R}(\phi) & = \sum_{\p} \Lambda^{(2), \ur}_{\q, \sigma_R}(\p) \widehat{\phi}\left(\frac{2 \log N\p}{\log c(\pi)}\right)\frac{2 \log N\p}{N\p\log c(\pi)} \\
& \qquad -  \frac{1}{2} \phi(0) + O\left(\frac{1}{\log c(\pi)}\right). \notag
\end{align}

\rk The extra contribution $\frac{1}{2} \phi(0)$ is crucial, and will be shown to be the only non-archimedean contribution to the type of symmetry. It is no surprise that the relations between Satake parameters and coefficients, as a specific property of each family, have an impact on the type of symmetry. This feature is already present in many classical works on low-lying zeros and its appearance when estimating the second moment has been explained by Dueñez and Miller \cite{duenez_effect_2009}.

\section{Trace formula}
\label{sec:stf}

\subsection{Selberg trace formula}
\label{subsec:fts}

Since the automorphic quotient of $G$ is compact, the original formulation of the trace formula\index{trace formula}, due to Selberg \cite{arthur_introduction_2005}, can be used and combined with the multiplicity one theorem. If $\Phi$ is a function in the Hecke algebra $\mathcal{H}(G(\A))$, then
\begin{equation}
\label{fts}
J_{\geom}(\Phi)= J_{\spec}({\Phi}),
\end{equation}

\noindent where the spectral and geometric parts\index{trace formula!spectral part}\index{trace formula!geometric part} are as follows. The geometric part is

\begin{equation}
\label{fts:geometrical-part}
J_{\geom}(\Phi) = \sum_{\{\gamma\}} \vol\left(G_\gamma(F) \backslash G_\gamma(\A)\right) \int_{G_\gamma(\A) \backslash G(\A)} \Phi\left(x^{-1}\gamma x\right) \dd x ,
\end{equation}

\noindent where the sum runs through conjugacy classes $\{\gamma\}$ in $G(F)$. Since $\Phi$ is compactly supported and $G(F)$ is discrete, the sum is finite. However its length depends on the support of $\Phi$ what turns to be a critical difficulty for estimations, for this support will depend on the spectral parameters, see Section \ref{subsec:ts-urpart}. The integrals appearing in this geometric side are called the orbital integrals\index{orbital integral, $\mathcal{O}_\gamma(\Phi)$}, defined by
\begin{equation}
\label{OI}
\mathcal{O}_\gamma(\Phi) = \int_{G_\gamma(\A) \backslash G(\A)} \Phi\left(x^{-1}\gamma x\right) \dd x .
\end{equation}

\noindent The spectral part is
\begin{equation}
\label{fts:spectral-part}
  J_{\spec}({\Phi})  = \sum_{\substack{\pi \subseteq L^2(G(F) \backslash G(\A))}} m(\pi) \widehat{\Phi}(\pi).
\end{equation}

\noindent Here $\pi$ runs through the isomorphism classes of unitary irreducible subrepresentations in $L^2(G(F) \backslash G(\A))$, and $\widehat{\Phi}$ is the Fourier transform of $\Phi$, see Section \ref{sec:plancherel-formulas}.

In order to have a problem amenable to the trace formula it is necessary to interpret statistics quantities on the universal family as a spectral side,  hence needed to select it by the Fourier transforms of suitable test functions. The aim of the present section is to construct a function $\Phi \in \mathcal{H}(G)$ such that, up to an error term,
\begin{align}
\label{selection}
J_{\mathrm{spec}}(\Phi) & \simeq \Lambda_{\d, \sigma_R}^{(\nu), \mathrm{ur}}.
\end{align}

\noindent In the case of factorizable test functions $\Phi = \otimes_v \Phi_v$, the spectral side of the trace formula factorizes as
\begin{equation}
\label{Fourier-transform:factorization}
\widehat{\Phi}(\pi) = \prod_v \widehat{\Phi}_v (\pi_v).
\end{equation}

\noindent Hence, in order to achieve the spectral selection \eqref{selection} it is sufficient locally select the conditions appearing in the decomposition of the universal family \eqref{uf-decomposition} through Fourier transforms. The following sections are dedicated to construct local test functions doing so, aim reached in Lemma \ref{lem:selecting}.

\subsection{Selecting the split conductor}
\label{subsec:split}

For an ideal $\d$ of $\mathcal{O}$, introduce the congruence subgroup given by the product of the corresponding local congruence subgroups in \eqref{filtration}, that is to say
\begin{equation}
\label{filtration-global}
K_0(\d) = \prod_{\p^r |\!| \d} K_{0, \p}(\p^r).
\end{equation}

The following result gives a test function whose Fourier transform selects the finite split conductor with multiplicities, in particular vanishes if $\c(\pi)$ does not divide a fixed ideal $\d$, see for instance \cite{lesesvre_counting_2020}.

\begin{lemma}
\label{lem:ft-finite-split-conductor}
For an ideal $\d$ of $\mathcal{O}$, let
\begin{equation}
\label{ft-finite-split-conductor}
\varepsilon_{\d} = \vol\left(\overline{K}_0(\d)\right)^{-1} \mathbf{1}_{\overline{K}_0(\d)}.
\end{equation}
Its Fourier transform selects the multiplicity relative to $\d$. More precisely, 
\begin{equation}
\label{Fourier-ft-finite-split-conductor}
\widehat{\varepsilon}_\d (\pi) = m(\pi, \d), \qquad \pi \in \mathcal{A}(G).
\end{equation}
\end{lemma}

\subsection{Selecting the ramified part}
\label{subsec:ramified}

For ramified places, the matrix coefficients allows to select the desired isomorphism class of supercuspidal representations in the following sense, see \cite[Corollary 10.26]{knightly_traces_2006}. Let $\sigma$ be a unitary representation of $G_R$. A matrix coefficient\index{matrix coefficient, $\xi_\sigma$} associated to $\sigma_R$ is a function of the form, given $v$ and $w$ in the space of $\sigma$, 
\begin{equation}
\label{matrix-coeff:def}
\begin{array}{cccc}
\displaystyle \xi_{\sigma}^{v, w} : &  G_R & \longrightarrow & \C \\
& g & \longmapsto & \langle \sigma(g)v, w \rangle.
\end{array}
\end{equation}

Matrix coefficients are continuous functions on $G_R$, compactly supported since $G_R$ is compact, and locally constant at finite places and smooth at archimedean places.

\rk The fact that matrix coefficients are considered only for ramified places is crucial for selecting purposes. The loss of the compactness of the support for matrix coefficients if there were archimidean split places, where non supercuspidal automorphic representations do exist, would make them fail to select the corresponding isomorphism class. Such a purpose can be achieved by means of existence theorems, for instance the Clozel-Delorme \cite{clozel_theoreme_1990} version of the Paley-Wiener theorem, at the expense of the precise form of the function. This is the reason why the non-totally definite case or the $\GL(2)$ case are analytically harder to deal with.

\begin{lemma}
Let $\sigma$ and $\pi$ be automorphic representations of $G_R$, and introduce $d_\pi$ the formal degree of $\pi$. Then for every unit vectors $v$ and $w$ in the representation space of $\sigma$, 
\begin{equation}
\label{matrix-coeff:orthog-rel}
\pi\left(\xi_{\sigma}^{v, w}\right) w = \mathbf{1}_{\pi \simeq \sigma} \frac{\langle w, v \rangle }{d_\pi} v.
\end{equation}
\end{lemma}

Taking for $v$ a vector of norm $d_\pi^{1/2}$, it follows that $\pi\left( \xi_\sigma^{v,v}\right)$ is the orthogonal projection onto $\C v$ and in the meanwhile selects the $\pi$'s isomorphic to $\sigma$. Considering its trace, this can be restated as follows.
\begin{lemma}
\label{prop:matrix-coeff}
Let $\sigma$ and $\pi$ be automorphic representations of $G_R$. Let~$v$ be a vector of norm one in the representation space of $\sigma$. Then,
\begin{equation}
\label{matrix-coeff:selecting-function}
\widehat{\xi_\sigma^{v,v}}(\pi) = \mathbf{1}_{\pi \simeq \sigma}.
\end{equation}
\end{lemma}

From now on, denote $\xi_\sigma$ any choice of matrix coefficient as in Proposition~\ref{prop:matrix-coeff}.

\subsection{The chosen test function}
\label{subsec:fonction-test}

Let $\d$ be an ideal of $\mathcal{O}^R$, $\sigma_R$ a representation of $\widehat{G}_R$ and $\p$ a prime ideal out of $R$. Introduce the test function
\begin{equation}
\Phi_{\d, \pi_R}^{(\nu), \ur}(\p) = \prod_v \Phi_v, 
\label{def:ft}
\end{equation}

\noindent which is built with the following local functions:
\begin{center}
\begin{tabular}{|L|L|}
\hline
\text{Places $v$} & \text{Local test function $\Phi_v$}\rule{0pt}{2.6ex} \\
\hline 
\notin R, \neq \p & \varepsilon_{\d, v}\rule{0pt}{2.6ex} \\
\hline
\in R, \neq \p & \xi_{\pi_v}\rule{0pt}{2.6ex} \\
\hline
\p & T_{\p^\nu}\rule{0pt}{2.6ex} \\
\hline
\end{tabular}
\end{center}

\noindent where
\begin{itemize}
\itex $\varepsilon_{\d}$ is the function introduced in Section \ref{subsec:split}, $\varepsilon_{\d, v}$ its $v$-component;
\itex $\xi_{\pi_v}$ is a matrix coefficient for $\pi_v$, see Section \ref{subsec:ramified};
\itex $T_{\p^\nu}$ is the Hecke operator, see Section \ref{subsec:hecke-opertaors-def}.
\end{itemize}

\begin{lemma}
\label{lem:selecting} Let $Q \geqslant 1$. Let $\d \wedge R = 1$ and $\sigma_R \in \widehat{G}_R$. Then
\begin{equation}
\label{Jspec, ftsapplied}
J_\spec\left(\Phi_{\d, \pi_R}^{(\nu), \ur}(\p)\right)  = \Lambda_{\q, \sigma_R}^{(\nu), \ur}(\p)+ O(N\p^{-\nu / 2}\Xi(\sigma_R)),
\end{equation}

\noindent where, introducing the set $X^\mathrm{ur}(G)$ of unramified characters of $G(\A)$,
\begin{equation}
 \Xi(\pi_R)  = \sum_{\substack{\chi \in X^\mathrm{ur}(G) \\ \chi_R \simeq \pi_R}} 1.
\end{equation}
\end{lemma}

\proof Let $\Phi = \Phi_{\d, \pi_R}^{(\nu), \ur}(\p)$. In order to determine the Fourier transform of~$\Phi$ recall that for every places $v$, $w$ and every $a \in \mathcal{H}(G_{v, w})$, $\widehat{a_va_w} = \widehat{a}_v\widehat{a_w}$. Thus,
\begin{equation}
\widehat{\Phi} = \prod_v \widehat{\Phi}_v = \widehat{T}_\p \prod_{\substack{v \in R \\ v \neq \p}} \widehat{\xi}_{\pi_v} \prod_{\substack{\mathfrak{f} \notin R \\ \mathfrak{f} \neq \p \\ \mathfrak{f}^r |\! | \d}} \widehat{\varepsilon}_{\f^r, v} .
\label{transform:ft}
\end{equation}

Hence only the Fourier transforms of the local components of the test function have to be determined.  The split part $\varepsilon_\d$ transforms into the characteristic function of conductors dividing $\d$ weighted by the corresponding multiplicities by Lemma~\ref{lem:ft-finite-split-conductor}. The ramified local parts $\xi_{\pi_v}$ transforms into the characteristic functions of the isomorphism class of $\pi_v$ by Lemma \ref{prop:matrix-coeff}. The Hecke operators transforms into the selecting function of unramified representation weighted by their associated Hecke eigenvalue by Lemma~\ref{lem:hecke}. The action of the Fourier transform of $\Phi$ follows, namely
\begin{equation}
\label{ft-Fourier-transformed}
\widehat{\Phi}(\sigma) =  m(\sigma, \d)  \lambda_\pi(\p^\nu) \1_{\substack{\sigma_R \simeq \pi_R \\ \p \nmid \c(\sigma^{R}) \div \d}}  .
\end{equation}

Nevertheless, these conditions also stand for characters: in order to not being killed by $\widehat{\Phi}$ they have to be trivial on $\overline{K}_0(\d)$, \textit{i.e.} they have to be unramified since $\det(\overline{K}_0(\d)) \subseteq \mathcal{O}^R$ and $\det(T(\p^\nu)) \subseteq \mathcal{O}_\p$. Moreover, they have to be isomorphic to $\pi_R$ at ramified places. The Fourier transform of the chosen test function hence does not vanish on unramified characters, unlike awaited. The corresponding extra contribution $\Xi(\pi_R)$ is treated separately in Lemma \ref{lem:characrers} below. After summing over the spectrum, it follows 
\begin{equation}
J_{\mathrm{spec}}(\Phi)  = \sum_{\substack{\sigma \in \mathcal{B}(\d, \pi_R)}} m(\sigma, \d) \lambda_\sigma(\p^\nu) + O\left(\sum_{\substack{\chi \in X^\mathrm{ur}(G) \\ \chi_R \simeq \pi_R}} 1 \right)
\end{equation}

\noindent and that achieves the proof. \qed

\subsection{Contribution of characters}

We now treat the extra contribution of characters. 
\begin{lemma}
\label{lem:characrers}
We have the following estimate, for every $\varepsilon > 0$ and  every~$Q \geqslant 1$,
\begin{equation*}
\frac{1}{Q^2} \sum_{\substack{N\q \leqslant Q \\ \q \wedge R = 1}} \sum_{\substack{\sigma_R \in \widehat{G}_R \\ c(\sigma_R) \leqslant Q/N\q}} \sum_{\d | \q} \lambda_2\left(\frac{\q}{\d} \right) \Xi(\sigma_R) \ll Q^{-1+\varepsilon}.
\end{equation*}
\end{lemma}

\proof Similarly to the intervention of the trace formula in Lemma \ref{lem:selecting} above, we use the Poisson summation formula in order to count characters. We seek a test function such that the spectral side gives $\Xi(\pi_R)$. We work on $F^\times$ instead of $Z \backslash F^\times$ for simplicity, characters of $Z \backslash F^\times$ corresponding to those of $F^\times$ trivial on the center. 

At archimedian places, the trivial action on the center implies they are among the trivial one and the sign, hence have conductor 1 at those places. Archimedean characters are of the form $\chi_{\varepsilon, t} = \mathrm{sgn}^{\varepsilon} |\det|^{it}$ for $\varepsilon \in \{0, 1\}$ and $t \in \R$.  Let us introduce $f_{T,v}$, for $v$ an archimedean place, a non-negative smooth function such that $\widehat{f}_{T,v}$ is compactly supported, takes values 1 for $t=0$, and $|\widehat{f}_{T,v}| \leqslant 1$. In particular, considering it as a function on archimedean characters, we get that $\widehat{f}_{T,v}(\mathrm{sgn}^{\varepsilon} |\det|^{it})$ is 1 for $t=0$, and vanish unless $t$ is small enough, say $|t| \leqslant T$.

At a finite split place $\p$, every character $\chi_\p$ can be written as $\chi_0 \circ \det$. Since the determinant of $\overline{K}_{0, \p}(\p^r)$ is $\mathcal{O}^{R \times}$, the effect of $\widehat{\varepsilon}_{\p^r}$ is to select characters $\chi_0$ unramified out of $R$. Introduce $f^R$ to be the characteristic function of ${\mathcal{O}^{R\times}}$, so that its Fourier transform selects unramified characters, by mimicking the proof of Lemma \ref{lem:ft-finite-split-conductor}. For ramified places, introduce $f_{\pi_v}$ to be the matrix coefficient associated to $\pi_v$, so that its Fourier transform selects characters locally isomorphic to $\pi_v$ as in Lemma \ref{prop:matrix-coeff}. Introduce 
\begin{equation}
f_{T, \pi_R} = f^R \prod_{v | \infty} f_{T,v}  \prod_{v \in R} f_{\pi_v} .
\end{equation}

Since the $\widehat{f}_{v}$ are non-negative and take value 1 on characters either unramified at split places or isomorphic to $\pi_v$ at ramified places, the Poisson formula gives
\begin{equation*}
\Xi (\pi_R) \leqslant \sum_{\chi \in \widehat{F^\times}} \widehat{f}_{T, \pi_R}(\chi) = \frac{1}{\vol(F^\times \backslash \A^\times)} \sum_{\gamma \in F^\times} f_{T, \pi_R}(\gamma).
\end{equation*}

Since $F^\times$ is a discrete set, choosing the archimedean components $f_{T, v}$ with a small enough support leads to kill every $f_{T, v}(\gamma)$ for $\gamma$ nontrivial. Hence $\Xi (\pi_R) \leqslant \vol(F^\times \backslash \A^\times)^{-1} f_{\pi_R}(1)$. The ramified part can be written as a Plancherel measure by the Plancherel inversion formula,
\begin{equation}
f_{\pi_R}(1) = \int_{\widehat{G}_R} \mathbf{1}_{\sigma \simeq \pi_R} \mathrm{d}\sigma = \mu^{\mathrm{Pl}}(\pi_R).
\end{equation}

\noindent Finally, coming back to the sum over all spectral data, we get 
\begin{align*}
& \frac{1}{Q^2} \sum_{\substack{N\q \leqslant Q \\ \q \wedge R = 1}} \sum_{\substack{\sigma_R \in \widehat{G}_R \\ c(\sigma_R) \leqslant Q/N\q}} \sum_{\d | \q} \lambda_2\left(\frac{\q}{\d} \right) \Xi(\sigma_R) \\
& \quad  \ll \frac{1}{Q^2} \sum_{\substack{N\q \leqslant Q \\ \q \wedge R = 1}} \sum_{\substack{\sigma_R \in \widehat{G}_R \\ c(\sigma_R) \leqslant Q/N\q}} \sum_{\d | \q} \lambda_2\left(\frac{\q}{\d} \right) \mu^{\mathrm{Pl}}(\sigma_R) \\
& \quad \ll \frac{1}{Q^2} \sum_{\substack{\sigma_R \in \widehat{G}_R \\ c(\sigma_R) \leqslant Q}}  \mu^{\mathrm{Pl}}(\sigma_R) \sum_{\substack{N\q \leqslant Q/c(\sigma_R) \\ \q \wedge R = 1}} \sum_{\d | \q} \lambda_2\left(\d \right) \\
& \quad \ll \frac{1}{Q^2} \sum_{\substack{\sigma_R \in \widehat{G}_R \\ c(\sigma_R) \leqslant Q}}  \mu^{\mathrm{Pl}}(\sigma_R) \sum_{\substack{N\q \leqslant Q/c(\sigma_R) \\ \q \wedge R = 1}} \mu(\q) \\ 
& \quad \ll \frac{1}{Q^{1-\varepsilon}} \sum_{\substack{\sigma_R \in \widehat{G}_R \\ c(\sigma_R) \leqslant Q}}  \frac{\mu^{\mathrm{Pl}}(\sigma_R)}{c(\sigma_R)^{1+\varepsilon}} 
\end{align*}

\noindent and this last sum over the ramified spectrum converges by \cite{lesesvre_counting_2020}. \qed

\section{Low orders contributions}
\label{sec:low-order}
 
\subsection{Unramified part}
\label{subsec:ts-urpart}

It remains to evaluate the sums of orders 1 and 2, displayed in \eqref{1-part} and \eqref{2-part}. Lemma \ref{lem:selecting} is the main tool for these evaluations, ultimately reducing the problem to the study of the geometric side of the trace formula.

\begin{prop}
\label{prop:ur-bound}
For $\phi$ an even Schwartz function whose Fourier transform has compact support in $(-2/3, 2/3)$, and for $\nu \in \{1, 2\}$,
\begin{equation}
\label{ur-bound}
\mathcal{P}^{(\nu), \ur}_Q(\phi) \ll \frac{1}{\log Q}.
\end{equation}
\end{prop}

The remainder of this section is dedicated to the proof of this proposition. The result obtained in Lemma \ref{lem:selecting} states that, up to an error term, the sought spectral sums $\Lambda_{\q, \sigma_R}^{(\nu), \mathrm{ur}}(\p)$ can be approximated by the spectral side \eqref{fts:spectral-part}. Recall that the Selberg trace formula states that this spectral part $J_\spec(\Phi_{\q, \pi_R}^{(\nu), \ur}(\p))$ is equal to the corresponding geometric side, which decomposes as
\begin{equation}
\label{Jgeom}
J_\geom\left(\Phi_{\q, \pi_R}^{(\nu), \ur}(\p)\right)  = J_1\left(\Phi_{\q, \pi_R}^{(\nu), \ur}(\p)\right) + J_\mathrm{ell} \left(\Phi_{\q, \pi_R}^{(\nu), \ur}(\p)\right),
\end{equation}

\noindent where the identity and elliptic terms are defined as
\begin{align*}
J_1\left(\Phi_{\q, \pi_R}^{(\nu), \ur}(\p)\right) & = \vol(G(F) \backslash G(\A)) \Phi_{\q, \pi_R}^{(\nu), \ur}(\p)(1) \\
J_\mathrm{ell} \left(\Phi_{\q, \pi_R}^{(\nu), \ur}(\p)\right) & = \sum_{\substack{\{\gamma\} \subset G(F) \\ \gamma \neq 1}} \vol(G_\gamma(F) \backslash G_\gamma(\A)) \int_{G_\gamma(\A) \backslash G(\A)} \Phi_{\q, \pi_R}^{(\nu), \ur}(\p) (x^{-1}\gamma x)\dd x .
\end{align*}

Since the identity lies outside the double classes $T(\p^\nu)$ defining the Hecke operator, $J_1(\Phi_{\q, \pi_R}^{(\nu), \ur}(\p))$ vanishes. The elliptic terms are bounded in the following lemma. 

\begin{lemma}
\label{lem:bounds-orbital}
For all $\gamma \in G(F)$, $\q$ ideal of $\mathcal{O}^R$, $\sigma_R \in \widehat{G}_R$ and $\nu \in \{1,2\}$,
\begin{equation}
\mathcal{O}_\gamma(\Phi) \ll_\varepsilon N\q^\varepsilon N\p^{\nu / 2} .
\end{equation}
\end{lemma}

\proof The orbital integrals factorize as a product of local ones: for a function $\Phi = \otimes_v \Phi_v$ in $\mathcal{H}(G)$, 
\begin{equation}
\mathcal{O}_\gamma(\Phi) = \prod_v \mathcal{O}_{\gamma_v} (\Phi_v).
\end{equation}

For ramified places, the associated local orbital integrals are uniformly bounded. This amounts to bound characters on $\mathrm{SU}(2)$ and this is a direct consequence of the Weyl character formula, see \cite{lesesvre_counting_2020}. For the finite split places, bounds are provided by Binder \cite{binder_fields_2017}: for a prime ideal $\mathfrak{f}$ different from $\p$ and every $\varepsilon > 0$,
\begin{equation}
\mathcal{O}_{\gamma_\mathfrak{f}}\left(\Phi^{\p, \nu}_{\q, \sigma_R, \mathfrak{f}}\right) \ll N\mathfrak{f}^\varepsilon .
\end{equation}

As for the $\p$-component of the test function, the orbital integral associated with the Hecke operator $T_{\p^v}$ is explicitly computed by Kottwitz \cite[Lemma 12.12]{matz_sato-tate_2016}, and yields
\begin{align*}
\mathcal{O}_{\gamma_\p}\left(\Phi^{\p, \nu}_{\q, \sigma_R, \p}\right) &= \mathcal{O}_{\gamma_\p}\left(T_{\p^\nu} \right) \\
& =  N\p^{-\nu /2} \mathcal{O}_{\gamma_\p}\left(\mathbf{1}_{T(\p^\nu)} \right)\\
& \ll N\p^{\nu/2}. \qed
\end{align*}

%


\noindent It is necessary to estimate the length of the sum over the elliptic conjugacy classes. In critical contrast with the previous article \cite{lesesvre_counting_2020}, for difficulties arise due to the presence of Hecke operators. Indeed, the support of the test function $\Phi$, in particular of the Hecke operators $T_{\p^\nu}$, are not uniformly supported in a compact as it is the case for the congruence subgroups. It is hence necessary to unveil the dependence on $\p$ when bounding the length of the sum and the global volumes. Even if general results due to Matz and Templier \cite{matz_sato-tate_2016} supply bounds in this case, the particular $GL(2)$-setting allows to be more precise and to obtain slightly better results. 
\begin{lemma}
\label{lem:contributing-classes}
The number of conjugacy classes $\{\gamma\}$ such that the orbital integral $\mathcal{O}_\gamma(\Phi_{\q, \pi_R}^{(\nu), \ur}(\p))$ is nonzero is uniformly bounded in $\q$ and $\sigma_R$. As for the $\p$-aspect, it is bounded by $N\p^{\nu / 2 + \varepsilon}$ for every $\varepsilon > 0$.
\end{lemma}

\proof The demonstration is a refinement of the counting argument provided by Matz \cite[Lemma 6.10]{matz_weyls_2017}. Let us lift the setting for convenience: an automorphic representation in $\mathcal{A}(G)$ is viewed as a cuspidal automorphic representation of $\PGL(2)$ by the Jacquet-Langlands correspondence, which is viewed as a cuspidal automorphic representation of $\GL(2)$ with trivial central character. 

Counting the $G(F)$-conjugacy classes is equivalent to counting the associated characteristic polynomials. Let $\gamma$ be a representative of a contributing conjugacy class, that is to say such that $\mathcal{O}_\gamma(\Phi)$ does not vanish. By definition of the test function, at all non-archimedean places $\gamma_v$ is a matrix with integer entries, for either $\gamma_\mathfrak{f} \in K_{0, \mathfrak{f}}(\mathfrak{r}^t)$ for a certain $t$ and an ideal $\mathfrak{r}$ prime to $\p$, or $\gamma_\p \in T(\p^\nu)$, both consisting of matrices with integer coefficients. Hence its characteristic polynomial $P_\gamma$ has coefficients in all the integers rings $\mathcal{O}_\mathfrak{r}$, hence in the integer ring of $F$. Introduce
\begin{equation}
P_\gamma = X^2 + a_\gamma X + b_\gamma, \qquad \alpha_\gamma, b_\gamma \in \mathcal{O}.
\end{equation}

Turning to the archimedean places, the test function is compactly supported modulo the center. Hence, up to normalizing the determinant to one  by replacing $\gamma$ by $\tilde{\gamma} = \gamma |\det \, \gamma|^{-1/2}$, the set of the contributing $\tilde{\gamma}$ lies in a fixed compact set, hence also the coefficients of the associated characteristic polynomials. In particular, the linear coefficient $a_{\tilde{\gamma}}$ is a bounded integer, and turning back to $\gamma$ we get
\begin{equation}
|a_\gamma |_\infty  \ll |\det \, \gamma|_\infty^{1/2}.
\end{equation}

The fact that $\gamma_\p$ lies in the Hecke double class $T(\p^\nu)$ and the others $\gamma_\q$ in the maximal compact subgroup $K_\q$ fixes the value of the non-archimedean norms of the determinant at each place , equal to $\p^\nu$ at the place $\p$, and to one at the other finite places. Since the determinant of $\gamma$ lies in $F$, the product formula yields
\begin{equation}
|a_\gamma|_\infty \ll |\det \gamma|_\infty^{1/2} = \prod_\mathfrak{f} |\det \, \gamma|_{\mathfrak{f}}^{-1/2} = |\det \, \gamma|_\p^{-1/2} = N\p^{\nu/2}.
\end{equation}

As for the factor $b_\gamma$, since it lies in the Hecke double class $T(\p^\nu)$, the same argument as above ensures that $b_\gamma$ is of archimedean norm equal to $N\p^\nu$. Therefore, by an application of Dirichlet's unit theorem, the number of choices for $b_\gamma$ is bounded by a power of $\log N\mathfrak{p}$.  \qed

Moreover, the global volumes have been precisely bounded by Matz \cite[Section 9]{matz_bounds} in the specific case of $GL(2)$: if $\gamma$ is in the support of $\Phi$, the volume $\mathrm{vol}(G_\gamma(F) \backslash G_\gamma(\A))$ is dominated by $N\p^\varepsilon$ for all $\varepsilon > 0$. The bounds on the number of contributing classes obtained in Lemma \ref{lem:contributing-classes} along with the bounds on orbital integrals at other places obtained in Lemma \ref{lem:bounds-orbital} imply that, for every $\varepsilon > 0$,
\begin{equation}
\label{lambda-contrib}
\Lambda_{\q, \sigma_R}^{(\nu), \ur}(\p) \ll N\p^{\nu+\varepsilon} \left(N\q\right)^{\varepsilon} + O \left( N\p^{\nu/2+\varepsilon} \Xi(\sigma_R) \right).
\end{equation}

The action of $\widehat{\phi}$ in the explicit formula, since compactly supported in $(-T_\phi, T_\phi)$, ensures a sum over primes running until $c(\pi)^{T_\phi/\nu}$ in the explicit formula \eqref{blo}. Plugging the above bounds in the definition of $\mathcal{P}^{(\nu), \ur}_{\q, \sigma_R}(\phi)$, the prime number theorem implies that
\begin{align*}
\mathcal{P}^{(\nu), \ur}_{\q, \sigma_R}(\phi) & \ll N\q^\varepsilon \sum_{\p}  |\widehat{\phi}|\left(\frac{\nu \log N\p}{\log c(\pi)}\right)\frac{\log N\p}{\log c(\pi)} N\p^{\nu/2+\varepsilon} \\
& \ll N\q^\varepsilon \frac{c(\pi)^{3 T_\phi / 2 + \varepsilon}}{\log c(\pi)}.
\end{align*}

After summation over $\q$, it is negligible compared to $|\mathcal{A}(Q)|$ if 
\begin{equation}
\frac{1}{Q^2}\sum_{\substack{N\q \leqslant Q \\ \q \wedge R = 1}} \sum_{\substack{\sigma_R \in \widehat{G}_R \\ c(\sigma_R) \leqslant Q/N\q}} \sum_{\d \div \q} N\d^\varepsilon \frac{c(\pi)^{3 T_\phi/2+\varepsilon}}{\log c(\pi)} \xrightarrow[Q\to\infty]{} 0,
\end{equation}

\noindent and this happens for $T_\phi \leqslant 2/3 - \varepsilon$ for any $\varepsilon > 0$, giving the desired result. \qed

%
%

\subsection{Ramified part}

It remains to estimate the contribution of $\p$-ramified representations to the spectral sum. This is the content of the following proposition.

\begin{prop}
\label{prop:r-bound}
For every $Q \geqslant 1$ and $\nu \geqslant 1$,
\begin{equation}
\label{r-bound}
\mathcal{P}_{Q}^{(\nu), \r}(\phi) \ll  \frac{Q^2}{\log(Q)^{\nu(1/2-7/64)}} \log \log Q.
\end{equation}
\end{prop}

\proof By definition of the conductor, $\pi$ is ramified at $\p$ if and only if $\p$ divides its arithmetic conductor. Hence, using the Blomer-Brumley bound and the counting law \eqref{cardinalityP},
\begin{align*}
\tilde{\Lambda}_{\q, \sigma_R}^{(\nu), \r}(\p) & = \mathbf{1}_{\p | \q} \sum_{\substack{\pi \in \mathcal{A}(\q, \sigma_R)}} \left( \alpha_\pi^\nu(\p) 
+ \beta_\pi^\nu(\p) \right) \\
& \ll \mathbf{1}_{\p | \q} N\p^{7\nu/64} A(\q, \sigma_R)   \\
& \ll  \mathbf{1}_{\p | \q} N\p^{7\nu/64}  \varphi_2(\q) \mu^\Pl(\sigma_R) + \mathbf{1}_{\p | \q} N\p^{7\nu/64}{R}(\q, \sigma_R).
\end{align*}

\noindent That leads to, after summing over the primes, 
\begin{equation}
\label{ram1}
\mathcal{P}^{(\nu), \r}_{\q, \sigma_R} (\phi) \ll \frac{ \varphi_2(\q) \mu^\Pl(\sigma_R)}{\log(N\q c(\sigma_R))}\sum_{\p | \q} \frac{\log N\p}{N\p^{\nu(1/2-7/64)}}.
\end{equation}

\begin{lemma}
\label{lem:sum-primes}
For every $0<s\leqslant 1$ and every $\q$,
\begin{equation}
\sum_{\p | \q} \frac{\log(N\p)}{N\p^s} \ll \log (N\q) ^{1-s} \log\log N\q.
\end{equation}
\end{lemma}
\proof This is a straightforward application of the hyperbola method. Indeed, for $Y>0$, partial summation gives
\begin{align*}
\sum_{\p | \q} \frac{\log(N\p)}{N\p^s} &= \sum_{\substack{\p | \q \\ N\p \leqslant Y}} \frac{\log(N\p)}{N\p^s} + \sum_{\substack{\p | \q \\ N\p > Y}} \frac{\log(N\p)}{N\p^s} \\ 
& \ll \sum_{N\p \leqslant Y} \frac{\log (N\p)}{N\p^s} + \sum_{\substack{\p | \q \\ N\p > Y}} \frac{\log(N\p)}{N\p^s} \\
& \ll \sum_{N\p \leqslant Y}N\p^{1-s}(\log (N\p+1) - \log N\p) + \frac{1}{Y^s} \sum_{\substack{\p | \q}} \log N\p \\
& \ll \max \left( Y^{1-s}, \frac{\log N\q}{Y^s} \right),
\end{align*}

\noindent and this quantity is optimized for $Y = \log N\q$, which gives the claimed statement. \qed

This applied to \eqref{ram1} imples that, for every $\q$ and $\sigma_R$,
\begin{equation}
\label{ram2}
\mathcal{P}^{(\nu), \r}_{\q, \sigma_R} (\phi) \ll  \frac{\varphi_2(\q) \mu^\Pl(\sigma_R)}{\log (N\q)^{\nu(1/2-7/64)}} \log \log N\q.
\end{equation} 

Applying now Lemma \ref{lem:dampening} when summing over the spectral data yields that for every $\varepsilon > 0$,
\begin{equation}
\label{r-bound}
\mathcal{P}_{Q}^{(\nu), \r}(\phi) \ll \frac{\log\log Q}{\log(Q)^{\nu(1/2-7-64)}}\sum_{\substack{\sigma_R \in \widehat{G}_R \\ c(\sigma_R) \leqslant Q}} \frac{\mu^\Pl(\sigma_R)}{c(\sigma_R)^{2-\varepsilon}},
\end{equation}

\noindent and this sum over the ramified spectrum converges by \cite{lesesvre_counting_2020}, ending the proof.  \qed

\noindent \textit{Remark.} If the automorphic representation $\pi$ is ramified at $\p$, then $\pi_\p$ is either a supercuspidal representation, a unitary twists of the Steinberg representation, or a unitary principal series. The two first cases are known to be tempered, so that the exponent $7/64$ could be replaced by $\varepsilon$. However, suitable unitary twists of a non-tempered unramified representation are ramified yet non-tempered. To rule out this case, a finer argument involving the support of the Fourier-Plancherel transform arising in the spectral side of the trace formula would be necessary. The above proof has the appeal of being both sufficient and uniform for all representations. We are grateful to the referee for having been the source of these comments.

In the explicit formula \eqref{ef}, there is a nontrivial contribution from the second order terms obtained in \eqref{2-part}, all the other terms being negligible by Propositions \ref{prop:ur-bound} and \ref{prop:r-bound}. Altogether this achieves the proof of Theorem~\ref{thmD}.

\section{Non-vanishing of $L$-functions}
\label{sec:non-vanishing}
The knowledge of the type of symmetry in Theorem \ref{thmD} leads to further statistics on the family of $L$-functions associated to representations in $\mathcal{A}(G)$. Following Iwaniec, Luo and Sarnak \cite{iwaniec_low_2000}, it opens the path to bounds on the density of non-vanishing at the central point. Introduce
\begin{equation}
\label{zero-order}
p_m(Q) = \frac{1}{|\mathcal{A}(Q)|} \# \left\{ \pi \in \mathcal{A}(Q) \ : \ \underset{s=1/2}{\mathrm{ord}} L(s, \pi) = m \right\}, \qquad m \geqslant 0.
\end{equation}

This section is dedicated to the proof of the following corollary. 
\nonvanishing*

\proof The proportion of vanishing at the central point could be estimated by approximating the Dirac mass $\phi = \delta_0$ in the one-level density limiting behavior \eqref{ts}. The Plancherel formula restates the asymptotic one-level density, for a Schwartz function $\phi$, as
\begin{equation}
\label{plancherel}
\int_\R \phi(x)W_O(x) \dd x = \int_\R \widehat{\phi}(y) \widehat{W}_O(y) \dd y.
\end{equation}

The proportion of vanishing at the central point, counted with multiplicities, can be bounded as follows. Introduce the functions
\begin{equation}
\phi(x) = \left(\frac{\sin \left(\pi T_\phi x\right)}{\pi T_\phi x} \right)^2 \qquad \text{and} \qquad \widehat{\phi}(y) = \frac{1}{T_\phi} \left( 1-\frac{|y|}{T_\phi} \right) \mathbf{1}_{|y| < T_\phi}
\end{equation}

\noindent as in \cite[equation (1.42)]{iwaniec_low_2000}.  Then $\phi$ is a non-negative function over $\R$ so that ${\phi}(0)=1$, \textit{i.e.} $\delta_0 \leqslant \phi$ over $\R$. The generalized Riemann hypothesis amounts to say that all the $\gamma_\pi$ are real so that, for all $Q \geqslant 1$,
\begin{align*}
\sum_{m \geqslant 1} mp_m(Q) & = \frac{1}{|\mathcal{A}(Q)|} \sum_{\pi  \in \mathcal{A}(Q)} \sum_{\gamma_\pi} \delta_0(\gamma_\pi) \\
& \leqslant  \frac{1}{|\mathcal{A}(Q)|} \sum_{\pi \in \mathcal{A}(Q)} \sum_{\gamma_\pi} \phi(\gamma_\pi) \\
& \leqslant  \frac{1}{|\mathcal{A}(Q)|} \sum_{\pi \in \mathcal{A}(Q)} D(\pi, \phi).
\end{align*}

\noindent Let $\phi$ be an even and Schwartz class function on $\R$, with Fourier transform supported in $(-2/3, 2/3)$. By the density result \eqref{ts} and the Plancherel formula \eqref{plancherel}, for every $\varepsilon > 0$ and for $Q$ sufficiently large, 
\begin{align*}
\frac{1}{|\mathcal{A}(Q)|} \sum_{\pi \in \mathcal{A}(Q)} D(\pi, \phi) 
& \leqslant \int_\R \phi(x)W_O(x) \dd x + \varepsilon \\
& \leqslant  \int_\R \widehat{\phi}(y) \widehat{W}_O(y) \dd y + \varepsilon.
\end{align*}

Iwaniec, Luo and Sarnak \cite[Appendix A]{iwaniec_low_2000} showed that the function $\phi$ chosen above is an optimal choice among functions supported in $(-T_\phi, T_\phi)$ for the orthogonal symmetry type. They computed 
\begin{equation}
\label{zero:values-g}
 \int_\R \widehat{\phi}(y) \widehat{W}_O(y) \dd y = \frac{1}{T_\phi} + \frac{1}{2},
\end{equation}

\noindent so that, since it is possible to chose a test function with Fourier transform compactly supported in $(-T_\phi, T_\phi)$ with $T_\phi = 2/3$ by Theorem \ref{thmD}, we deduce from the above,  by letting $\varepsilon$ go to zero, that
\begin{equation}
\liminf_{Q \to \infty} \sum_{m \geqslant 1} mp_m(Q) \leqslant \frac{1}{2} + \frac{1}{T_\phi} = 2,
\end{equation}

\noindent and this proves Corollary \ref{coro:coro}. \qed

\noindent \textit{Remark.} The very definition of these proportions gives the for every $Q\geqslant 1$,
\begin{equation}
\label{zero:sum-pm}
\sum_{m \geqslant 0} p_m(Q) = 1.
\end{equation}

\noindent This implies a lower bound for the non-vanishing proportion
\begin{align*}
p_0(Q) & \geqslant \sum_{m \geqslant 0} p_m(Q) - \sum_{m \geqslant 0} mp_m(Q) \\
& \geqslant 1 - \int_\R \widehat{\phi}(y) \widehat{W}_O(y) \dd y - \varepsilon,
\end{align*}

\noindent so that using the optimal choice of function $\phi$ as above it yields, for $\widehat{\phi}$ supported in $(-T_\phi, T_\phi)$,
\begin{equation}
\liminf_{Q \to \infty} p_0(Q)  \geqslant \frac{1}{2} - \frac{1}{T_\phi}, \\
\end{equation}

\noindent providing a family of bounds depending on the function $\phi$. The bound on the support of the Fourier transform of $\widehat{\phi}$ by $2/3$ is too small to yield a nontrivial result on $p_0(Q)$. Enlarging the support for the Fourier transform is expected to be a challenge since off-diagonal terms enter into account in the trace formula. 

\bibliographystyle{cdraifplain}
\bibliography{biblio}
\end{document}